\documentclass[12pt]{article}
\input{epsf.sty}
\usepackage{hyperref}
\usepackage{amsmath, amssymb}

\usepackage[arrow, matrix, curve]{xy}

\newtheorem {thm}{Theorem}[section]
\newtheorem {prop}[thm]{Proposition} 

\newtheorem {lem}[thm]{Lemma}
\newtheorem {cor}[thm]{Corollary}

\newcommand{\qed}{\nobreak \ifvmode \relax \else
      \ifdim\lastskip<1.5em \hskip-\lastskip
      \hskip1.5em plus0em minus0.5em \fi \nobreak
      \vrule height0.75em width0.5em depth0.25em\fi}

\def\Cox{\hfill \Box}

\def\N{{\Bbb N}}

\def\Z{{\Bbb Z}}
\def\R{{\Bbb R}}
\def\LL{{\Bbb L}}

\def\E{{\Bbb E}}

\def\e{{\varepsilon}}

\def\ba{{\backslash}}
\def\sb{{\subset}}

\def\D{\Delta}
\def\a{\alpha}

\def\ba{\setminus}
\def\b{\beta}

\def\d{\delta}

\def\e{\varepsilon}

\def\phi{\varphi}

\def\g{\gamma}

\def\l{\lambda}

\def\LL{\mathcal{L}}
\def\s{\sigma}

\def\t{\tau}

\def\th{\theta}

\def\o{\omega}

\def\D{\Delta}

\def\L{\Lambda}

\def\O{{\Omega}}

\def\T{\T}

\def\HH{{\cal H}}

\def\GG{{\cal G}}
\def\PP{{\cal P}}

\def\II{{\cal I}}

\def\V|{{\Vert}}

\begin{document}
\title{A class of non-ergodic probabilistic cellular automata with unique invariant measure and quasi-periodic orbit}
\author{
Benedikt Jahnel
\footnote{ Ruhr-Universit\"at   Bochum, Fakult\"at f\"ur Mathematik, D44801 Bochum, Germany,
\newline
 \texttt{Benedikt.Jahnel@ruhr-uni-bochum.de}, 
\newline
\texttt{http://http://www.ruhr-uni-bochum.de/ffm/Lehrstuehle/Kuelske/jahnel.html }}
 \, and  Christof K\"ulske
\footnote{ Ruhr-Universit\"at   Bochum, Fakult\"at f\"ur Mathematik, D44801 Bochum, Germany,
\newline
\texttt{Christof.Kuelske@ruhr-uni-bochum.de}, 
\newline
\texttt{http://www.ruhr-uni-bochum.de/ffm/Lehrstuehle/Kuelske/kuelske.html
/$\sim$kuelske/ }}\, 
\,  
}

\maketitle

\begin{abstract} 
{\bf }
We provide an example of a discrete-time Markov process on the three-dimensional infinite integer lattice with 
$\Z_q$-invariant Bernoulli-increments which has as local state space the cyclic group $\Z_q$. 
We show that the system has a unique invariant measure, but  remarkably 
possesses an invariant set of measures on which the dynamics is conjugate to an irrational rotation on the continuous sphere $S^1$. 
The update mechanism we construct is exponentially well localized on the lattice.

\end{abstract}



\vfill\eject

\vfill\eject

\newpage 
\section{Introduction}

The possible long-time behavior of infinite lattice systems under stochastic dynamics is a subject of ongoing research. 
The situation is non-trivial already for reversible dynamics and 
becomes even more difficult when one leaves the assumption of reversibility 
of the dynamics and enters the realm of driven systems. 
Infinite lattice systems
%, unlike finite-state \textbf{Markov chains},  
may possess different equilibria for the same stochastic dynamics. 
  The first question which comes to mind is to estimate 
the approach to equilibrium, if there is a unique equilibrium (using for example spectral gap analysis, logarithmic Sobolev inequalities, etc. see \cite{Ch04,Ch05}).   
If there are multiple equilibria one 
may be interested in their domains of attraction.  
An interesting question in this context is whether a unique equilibrium has to be an attractor for a stochastic lattice dynamics (see for example \cite{Li85,ChMa11,Mo95,RaVa96}).  
Ideally one would like to understand possible behavior  of invariant sets and attractors. 
Under what circumstances can there be oscillatory behavior and closed orbits of measures?  
This is difficult to answer for infinite lattice systems.  

On the one hand motivation comes from the study of given models describing e.g systems of coupled neurons. These show characteristic patterns of spiking, phases of long-range order, periodicity, synchronization as well as disordered phases. A mathematical analysis 
of these interesting phenomena so far has been restricted mainly to mean-field models (see for example \cite{DeGaLoPr14} or 
%phenomenon of  of systems of interacting oscillators, proved to occur in a 
the Kuramoto model \cite{DH96,AcBoPeRiSp05,GiPaPePo12,GPP12}). 
On the other hand there is the theoretical interest to make progress by restricting the possible 
forms of limiting behavior, or to explore possible forms and 
provide model systems which illustrate possible forms of "non-standard" limit behavior. 
Ultimately one may strive for a classification of models into classes which look similar, but in the beginning we are forced 
to work with model systems, which serve as illustrations, their universality being a question that hopefully can be tackled later.  

\medskip
These types of question may be posed for continuous-time Markov processes or discrete-time Markov processes. 
In a previous paper \cite{JaKu12} we were able to settle an old question raised by Liggett \cite{Li85} whether it is possible to have a continuous-time Markovian dynamics with
a unique invariant measure which does not attract all initial conditions, but has a closed 
orbit of measures. 
An example of a probabilistic cellular automaton (PCA) which shows that non-ergodicity is possible in one dimension with positive rates was given in \cite{Ga01}. 
An example of a PCA which shows that non-ergodicity in one dimension, even with a unique stationary measure, is possible was given in \cite{ChMa11} (albeit with some deterministic updatings).   
%%A PCA is a Markov process on the state space $E^{\Z^d}$ where $E$ is a finite set. 
%They have $\{0,1\}$ as the local state space and devised their dynamics in such a way that the two perfect alternating configurations $\o_{01}=\dots, 0,1,0,1,0 \dots$ (by which we denote the alternating configuration of zeros and ones with a zero at the origin), and its shift by one lattice site, $\o_{10}$ reproduce each other under iteration of the dynamics (they "rotate"). Configurations which differ locally from the alternating ones undergo a random transition mechanism which tends to bring them closer to the alternating ones. 
%Further they prove 
%%Chassaing and Mairesse then prove in their paper 
%that the only stationary measure indeed is the symmetric mixture $\frac{1}{2}(\d_{\o_{10}}+\d_{\o_{01}})$ of the Dirac measures. Since there is a periodic orbit of length two given by these two Dirac measures, the system is not ergodic. It would be nice nevertheless to have a less degenerate example which does not need that some transitions are deterministic, as above. Our dynamics has the property of being purely non-deterministic although it is only a weak PCA with quasilocal updating. 
In \cite{DaLoRo02} another example of a non-ergodic, non-degenerate, Ising-type PCA on the $d$-dimensional lattice with $d\geq2$ is exhibited.
%$\{-1,1\}^{\Z^d}$ is exhibited. 
%Similarly to our dynamics in the phase-transition region two Gibbs measures $\nu^+\neq\nu^-$ alternate under the PCA. Further $\a\nu^++(1-\a)\nu^-$ is time-stationary if and only if $\a=1/2$. They do not provide any reasoning that $\frac{1}{2}(\nu^++\nu^-)$ is maybe the unique time-stationary measure for their model.
Both examples have periodic orbits of period two.
To compare, in the case of an interacting particle system (IPS) where the updating is in continuous time, non-ergodicity with unique invariant measure was proved to be impossible for local rates in one lattice dimension in \cite{Mo95,RaVa96}. In two lattice dimensions the question remains open. For general background, see also \cite{LeMaSp90,To01}.

\medskip
The aim of the present paper is to provide an example of a discrete-time Markov process with 
{\em discrete} $\Z_q$-invariant Bernoulli-updates having a unique invariant measure,  
but a non-trivial quasi-periodic orbit under the dynamics. 
It is remarkable that a quasi-periodic orbit which is conjugate to an irrational rotation on a sphere can arise from exponentially localized dynamical rules on a discrete local spin space.

To do so, we consider a class of discrete $q$-state spin equilibrium models defined in terms of a  translation-invariant 
quasilocal specification with discrete clock-rotation invariance having extremal Gibbs measures $\mu'_{\phi}$ labeled by the uncountably many values of $\phi$ in the one-dimensional sphere (introduced 
%by  van Enter, Opoku and K\"ulske 
in \cite{EnKuOp11}).
Next we 
construct an associated discrete-time Markov process as a 
PCA with exponentially localized updating rule. 
The process has the property to reproduce a deterministic rotation of the extremal Gibbs measures and preserve macroscopic coherence by a quasilocal and time-synchronous updating mechanism in discrete time without deterministic transitions.
% (microscopic random updating).
We prove that,  depending on an updating-velocity parameter $\t$ there is either a continuum of time-stationary measures and closed orbits of rotating states, or a unique time-stationary measure and a dense orbit of rotating states. In both cases the process in non-ergodic. 

Our paper partially builds on constructions of \cite{JaKu12} where we considered an IPS which can be seen as a  
continuous-time analogue of the present Markov process.   
%and proved that the former has a unique translation-invariant time-invariant measure, but possesses non-trivial closed orbits of measures given by rotating states. 
The bulk of the present work consists in the proof of uniqueness of the invariant measure in the absence of continuous time, 
which heavily draws on techniques relating entropy loss, the inverse transition operator and the Gibbs property of the invariant measure.

\subsection{Discrete-time infinite lattice models -- Main result}

%$\cal h$

We consider Markovian dynamics on the infinite-volume state space 
$\O=E^G$, where $E$ is a finite local state space sitting at each site $i\in G$ and $G$ is a countable set.  
In the examples we want to construct we will specifically choose 
$E=\Z_q=\{0,1,\dots, q-1\}$ to be the cyclic group and 
$G=\Z^d$ to be the $d$-dimensional integer lattice, so that the infinite volume configuration space
$\O$ is an abelian group 
w.r.t. sitewise addition modulo $q$. 

Recall the following definitions.
A {\em (deterministic) cellular automaton} is given by a deterministic local updating rule  
$f: E^\L \mapsto E$ where $\L\subset G$ is a finite set containing the origin $\L\ni 0$. 
The simultaneous application of $f$  yields a corresponding function $F:\O\mapsto \O$ acting on infinite-volume configurations $\s\in \O$, where the $i$-th coordinate of the image $F(\s)$ is 
given by $F(\s)_i=f(\th_i\s)$ for all sites $i \in G$ 
and $\th_i \s$ is the configuration obtained from $\s$ by a lattice shift by the vector $i\in G$, that is  
$(\th_i \s)_j = \s_{j+i}$. One then likes to study properties of the discrete dynamical system obtained by 
iterates of $F$. 
We will not discuss such deterministic cellular automata any further, but mention that  
there are recent developments concerning very interesting questions about simple cellular automata which are fed 
random initial configurations, bootstrap percolation being 
one of these (see for example \cite{AiLe88,BoDuMoSm14,BaBoMo09,DuEn13}).

A {\em (strict) probabilistic cellular automaton} (PCA) is 
given by a probabilistic single-site updating kernel $p$ from the infinite-volume 
state space $\O$ to single-site probability measures,  
which plays the role of a random version of the map $f$. 
So one has  
$p: \O \mapsto \PP(E)$ where $\PP(E)$ denotes the probability measures on the finite set $E$
and $p$ is assumed to be a strictly local probability kernel, which means that for each $a\in E$ the functions 
$\s\mapsto p(\s,a)$ are strictly local functions for all $a\in E$, that is functions depending 
only on finitely many coordinates of $\s$. 
This probabilistic updating rule is applied (stochastically) independently over all sites 
to yield a Markovian transition operator $M$ on infinite lattice 
configurations (which generalizes the deterministic function $F$) of the form 
\begin{equation*}\label{Product}
 \begin{split}
M(\s,d\eta)=\otimes_{i\in G}p(\theta_i\s, d\eta_i).
\end{split}
\end{equation*}
This transition operator describes 
simultaneous updates on the infinite lattice which are performed in a sitewise independent way according 
to a lattice-shift invariant rule. 
There is a huge literature on probabilistic cellular automata (see for example \cite{ChMa11,LeMaSp90,Ga01,DaLoRo02,To01,ToVaStMiKuPi90}).

A {\em probabilistic cellular automaton with exponentially localized update kernel} 
is given by a Markovian kernel $M$ from the set of infinite-volume configurations 
$\O$ to the probability measures on $\O$ 
%with the following two properties.  
%They 
allowing for exponentially suppressed non-localities in the following sense. 

1. The update kernel $M(\s, d\eta)=:M_\s(d\eta)$ is {\em uniformly spatially mixing in the future with exponent
 $\a_2$}. 
  This means by definition that there is a positive 
 constant $K_2$, such that  at any fixed starting configuration $\s$ 
 the following decay estimate holds:  
For all finite volumes $\L \subset \D \subset G$,  for all sets $A$ in the sigma-algebra 
%measurable w.r.t. the sigma-algebra 
generated 
by the $\eta$-coordinates in $\L$ and sets $B$ in the sigma-algebra generated by the $\eta$-coordinates in $G\setminus\D$ one has 
\begin{equation}\label{in the future}
 \begin{split}
&|M_\s(A|B) - M_\s(A)| \leq K_2  \sum_{i\in \L, j\in G\setminus\D }e^{- \a_2 |i-j| }.
\end{split}
\end{equation}
We remark that this uniform mixing property is much stronger (closer to independence) than simply correlation decay.  
While correlation decay 
holds e.g whenever the conditional measure on the $\eta$ variables is an extremal Gibbs measure even 
in the phase-transition regime, phase transitions are excluded by \eqref{in the future}.

2. The update kernel $M$ satisfies the property of {\em exponential locality from the past with exponent $\a_1$}. 
This means by definition that there is a positive constant $K_1$ such that 
the variation w.r.t. $\s_i$ of the probability $M(\s,\eta_\L)$ with finite $\L$ 
%to find a configuration $\eta_V\in E^V$ in the finite volume $V$ in an update from the infinite-volume configuration $\o$ 
is bounded by the exponential estimate 
\begin{equation*}\label{Product}
 \begin{split}
&\d_{i}( M(\cdot ,\eta_\L)):= \sup_{{\s,\tilde \s:}\atop{\s_{i^c}=\tilde\s_{i^c}}}(M(\s,\eta_\L)-M(\tilde \s,\eta_\L))
%\leq K_1\sum_{j\in\L} e^{-\a_1|i-j|}
\leq K_1e^{-\a_1|i-\L|}.
\end{split}
\end{equation*}
Here we used notation $i^c=G\setminus \{i\}$, $|i-\L|=\min_{j\in \L}|i-j|$ and $\eta_\L:=\{\xi\in\O: 1_{\eta_\L}(\xi)=1\}$.

Notice that kernels with these two properties are related to the infinite-volume Markovian kernels studied in \cite{Ku84}.

\medskip
Specifying now to Markov kernels $M$ on the state space $(\Z_q)^G$ we say that $M$ 
is {\em $\Z_q$-invariant} iff it is compatible with 
joint rotation on the local state space and we can write $M(\s+a, A+a)= M(\s, A)$ for all $a \in \Z_q$. 
We say that $M$ has {\em Bernoulli-increments} if $M(\s, \s+\{0,1\}^G)=1$. 
In that case the updated configuration is obtained from an initial configuration $\s$ 
by the site-wise addition modulo $q$ of a $\{0,1\}^G$-valued 
field $N$ of Bernoulli-increments whose distribution conditional on $\s$ we denote 
by the symbol with a hat, namely $\hat M(\s, dn)$. 
To exclude degeneracies we only consider kernels in this paper which are {\em uniformly non-null}  meaning by definition that 
$\hat M_\s(n_{\L})\geq c^{|\L|}$ for all $n_{\L}\in \{0,1\}^{\L}$ 
for some strictly positive uniform constant $c$.

Given such a Bernoulli updating kernel  
$\hat M(\s, dn)$ from $(\Z_q)^{\Z^3}$ to $\{0,1\}^{\Z^3}$,  
we then look at 
the associated discrete-time Markov process 
\begin{equation*}\label{Product}
 \begin{split}X^t=(X^{t-1}+N^{t})\text{ mod } q
 \end{split}
\end{equation*}
 on the state space $\O=(\Z_q)^{\Z^3}$ 
with increment distribution 
\begin{equation*}\label{Product}
 \begin{split}\LL\bigl(N^t=\cdot \bigl |X^{t-1}=\s\bigr)=\hat M(\s, \cdot)
 \end{split}
\end{equation*}
%{\bf HOW CAN WE WRITE LAW IN A BETTER FONT?}
and use the symbol $M$ (without the hat) for the corresponding transition kernel, i.e.  
\begin{equation*}\label{Product}
 \begin{split}M(\s, \cdot )=\LL\bigl(X^t=\cdot \bigl |X^{t-1}=\s\bigr).
 \end{split}
\end{equation*}
We write $\PP_{\theta}(\O)$ for the lattice-translation invariant probability measures on $\O$. 

\medskip
The main result of the paper is the following theorem which states that we may 
construct exponentially localized transition kernels with Bernoulli updates
which are arbitrarily well localized, 
but whose time-evolutions show macroscopic coherence and non-ergodicity 
in the PCA-sense. 

\begin{thm}\label{PCA_Theorem}  For any arbitrarily large prescribed mixing exponents 
$\a_1,\a_2\in (0,\infty)$ there exists  
an integer $q_0\in \N$ such that, for all $q\geq q_0$ 
the following is true. 

There exists an updating kernel $M$ with Bernoulli increments  
on $(\Z_q)^{\Z^3}$,  
which satisfies the properties of uniform exponential mixing in the future with an exponent of (at least) 
$\a_2$ and uniform spatial locality from the past with an exponent of (at least) $\a_1$, 
such that the associated discrete-time stochastic dynamics 
possesses a set $\II\subset \PP_{\theta}(\O)$ of lattice translation-invariant measures 
%which is invariant for the dynamics 
on which the dynamics acts like a rotation. More precisely: 
\begin{enumerate}
\item
There is a quasilocal observable $\psi:\O\rightarrow \R^2$ such that the expectation map 
$$\nu \mapsto \nu(\psi)$$ is a bijection from $\II$ to the sphere $S^1\subset \R^2$.    
\item The stochastic dynamics restricted to the invariant set 
$\II$ is conjugate to a rotation by some angle $\t=\t(M)$: 
$$\forall \nu \in \II:\quad (M\nu)(\psi)=\nu(\psi)+ \t$$ 
where we have written rotation on $S^1$ on the r.h.s in an additive way.

\item
The  kernel $M$ can always be chosen so that the cases, $\t/2\pi$ is rational, or irrational can both occur. 
\begin{enumerate}
\item  If $\t/2\pi$ is irrational, there exists a unique invariant measure $\nu^*$ among the translation-invariant 
measures, 
but the dynamics is quasiperiodic. In particular $M^n \nu \not\rightarrow \nu^*$ (in the sense of local convergence)
and the PCA is non-ergodic.
%in the {\bf PCA-sense}. 
  
\item If $\t/2\pi$ is rational, there are uncountably many periodic orbits and uncountably 
many invariant measures. 
\end{enumerate}

\end{enumerate}

\end{thm}

\subsection{Comparison with continuous-time dynamics}\label{Comp}

We believe that it is instructive to compare this result to our previous work 
on the existence of non-ergodic  
{\em continuous-time} dynamics in \cite{JaKu12}. For easy comparison let us also present the following result which is based on the work of \cite{JaKu12}.

We consider a continuous-time interacting particle system (IPS) on $(\Z_q)^{\Z^3}$ given in terms 
of a generator $U$ acting on observables $\psi:(\Z_q)^{\Z^3}\mapsto \R$
of the form  
\begin{equation}\label{IPS}
U\psi(\o)=\sum_{i \in \Z^d}\Bigl( c_i^+(\o)\bigl(\psi(\o+1_i)-\psi(\o) \bigr)  +c_i^-(\o)\bigl(\psi(\o-1_i)-\psi(\o) \bigr)     \Bigr)
\end{equation}
where $(\o\pm 1_i)_j=(\o_j \pm  \d_{i,j})$ mod $q$, so that it allows both for increments $\pm 1$ at each site $i$.  

In analogy to the notion of uniform exponential locality from the past in the discrete-time setup 
we say that the rates $ c_i^\pm (\o)$ satisfy the property of {\em exponential locality with exponent $\a_1$} whenever 
there is a finite $K_1$ such that 
\begin{equation}\label{Exp_Loc}
 \begin{split}
&\d_{i} c_j^\pm(\cdot):
= \sup_{{\o,\tilde \o:}\atop{\o_{i^c}=\tilde\o_{i^c}}}
(c_j^\pm(\o)- c_j^\pm(\tilde\o) )\leq K_1 e^{- \a_1 |i - j |}.
\end{split}
\end{equation}
This in particular makes the dynamics well-defined by standard methods. 
%Corresponding to the notion of uniform spatial mixing in the future we have the following weaker property of \textit{exponential decay of correlations with exponent $\a_2$}
%\begin{equation}\label{Dec_Cor}
% \begin{split}
%\sup_{\s}|S_\t1_{\eta_{\L\cup\bar\L}}(\s) - S_\t1_{\eta_{\L}}(\s)S_\t1_{\eta_{\bar\L}}(\s)| 
%\leq K_\t\sum_{i\in\L,j\in\bar\L}e^{-\a_2|i-j|}
%\end{split}
%\end{equation}
%where $S_\t$ is the Markov semigroup associated to the generator $L$, $\t>0$ and $\L,\bar\L\subset G$ finite volumes.
Then we have the following theorem. 

\begin{thm} \label{PCA_IPS}  For arbitrarily large exponent
$\a_1\in (0,\infty)$ there exists  
an integer $q_0\in \N$ such that, for all $q\geq q_0$ 
the following is true: There exists a generator
$Q$ on $\O=(\Z_q)^{\Z^3}$ with lattice-translation invariant and $\Z_q$-invariant rates which 
satisfy exponential locality with exponent (at least) $\a_1$
% and exponential decay of correlations with exponent (at least) $\a_2$
such that the associated Markov process with continuous-time semigroup $(S_t)_{t\geq0}$ 
possesses a set $\II\subset \PP_{\theta}(\O)$ of lattice translation-invariant measures
%which is invariant for the dynamics 
on which the dynamics acts like a rotation. More precisely: 

\begin{enumerate}
\item (This is identical to the corresponding point for the discrete-time dynamics.) 
There is a quasilocal observable $\psi:\O\rightarrow \R^2$ such that the expectation map 
$$\nu \mapsto \nu(\psi)$$ is a measurable bijection from $\II$ to the sphere $S^1\subset \R^2$.    
\item The continuous-time stochastic dynamics restricted to the invariant set 
$\II$ is conjugate to a continuous rotation by $t\in \R$:  
$$\forall \nu \in \II:\quad (S_t\nu)(\psi)=\nu(\psi)+ t.$$ 
\item There is a unique time-stationary measure $\nu^*$ among the lattice-translation invariant measures, 
namely the uniform mixture over $\II$. 
 \end{enumerate}

\end{thm}

%\begin{rem}
Note that in the case of the continuous-time Markov process we do not need a requirement analogous to property \eqref{in the future} posed in the PCA case. 
%Nevertheless for $S_t$ with $t\geq0$ we do have a similar although weaker condition on the covariances see Proposition \ref{CorrelationDecay}.

Neither does Theorem \ref{PCA_Theorem} imply Theorem \ref{PCA_IPS} nor does Theorem \ref{PCA_IPS} imply Theorem \ref{PCA_Theorem}.

\subsection{Ideas of the proof}

There are common parts (see A below) and essential differences 
(see C, and the more difficult D) in the treatment of discrete-time dynamics 
and continuous-time dynamics.  

A. First we give a one-parameter family of measures $\II= \II(\b,q)\subset \PP_{\theta}((\Z_q)^{\Z^3})$, 
depending on an inverse temperature parameter $\b$ for each sufficiently large $q$. 
This family is used in both cases of discrete-time and continuous-time dynamics. 
We will shortly review its construction which was given in \cite{JaKu12}, based on arguments concerning 
the preservation of Gibbsianness, in Section \ref{The equilibrium model}.  

B. We construct a discrete-time Bernoulli update kernel for which  $\II$ is an invariant set in Section \ref{The updating mechanism}. 
This is analogous to but different from the construction of $U=U(\b,q)$ given in \cite{JaKu12}. 

C. We prove locality, mixing, and further properties also in Section \ref{The updating mechanism}. 

D. We prove uniqueness of the invariant measure 
(where  it is claimed to hold). The discrete case does not follow 
from the continuous case since time-derivatives are not available, and it necessitates the use of a different chain of arguments. 

\medskip

Physically the construction of the rotating-states mechanism
%of a
%this non-ergodic 
%discrete-time Markov process 
is inspired by conjectures in \cite{MaSh11} in the context of  IPS based on a clock model in an intermediate-temperature regime \cite{FrSp82}. 
%It is also related to the phenomenon of synchronization of systems of interacting oscillators, proved to occur in a mean-field setup for the Kuramoto model, see \cite{DH96,AcBoPeRiSp05,GiPaPePo12,GPP12}. 
%The rotating-states mechanism also was the inspiration for our related but different 
%construction \cite{JaKu12,JaKu14} of a non-ergodic IPS with unique stationary measure. 
To carry out our construction, as in \cite{JaKu12,JaKu14} we draw on the relation to the planar rotator model which has as a local state space the one-dimension sphere $S^1$. On the lattice in three or more space dimensions, at sufficiently strong coupling constant this system exhibits the breaking of the rotation symmetry in spin-space, see \cite{FrSiSp76,Pf82,FrPf83,MaSh11}. To arrive at a system of discrete spins (or particles) with finite local state space 
%$\{0,\dots, q-1 \}$
$\Z_q$ a local discretization is applied for $q$ sufficiently large but finite. Then the interplay between the systems of discrete and continuous spins is exploited. In particular we use the fact that the discretization map bijectively maps the lattice translation-invariant extremal Gibbs measures $\mu_{\phi}$ of the continuous system to the extremal lattice translation-invariant Gibbs measures $\mu'_{\phi}$ of the discrete system where  $\phi$ runs over the one-dimensional sphere $S^1$. Note also the non-trivial fact that the discrete system has uncountably many extremal Gibbs measures.   

As we will see below we can define an associated discrete-time Markov process with transition kernel $M_\t$, where $\t$ is a continuous parameter carrying also the meaning of an angle. 
The kernel $M_\t$  assigns to  a particle configuration 
in $(\Z_q)^{\Z^d}$ a random particle configuration in $(\Z_q)^{\Z^d}$ and will be obtained by a natural three-step procedure. 
We call this procedure {\em Sample-Rotate-Project}, to be described in detail below. Here it is the rotation step which 
carries the dependence on the continuous angle $\t$. For $0\leq\t\leq2\pi/q$ the updating is Bernoulli.
%The present construction of a discrete-time Markovian dynamics also helps to understand our previous construction of an non-ergodic IPS, as we will explain, since the latter can be seen as an infinitesimal-$\t$ version of the former. 
The dynamics works nicely on the Gibbs measures and we have the rotation property: 
\begin{enumerate}
\item An application of the transition operator $M_\t$ to a discrete Gibbs measure $\mu'_{\phi}$ 
yields a rotation by an angle $\t$, so that we have $M_{\t}\mu'_{\phi}=\mu'_{\phi+\t}$. 
\end{enumerate}
The proof of this fact is more straightforward than the proof of the analogous statement in the IPS setup \cite{JaKu12}, given the previous work on preservation of Gibbsianness under discretizations. Property 1 already implies that the symmetric mixture $\mu'_*=\frac{1}{2\pi}\int_0^{2\pi} d\phi\mu'_{\phi}$ is invariant under the dynamics. 

Note also that we can play with the velocity-parameter $\t$ now, and consider the action of the dynamics on the Gibbs measures.  
Rational values of $\t/2\pi$ yield finite closed orbits, of which there are uncountably many, 
so that there are uncountably many time-stationary measures obtained as the equal-weight measures on these orbits. 
Irrational values of $\t/2\pi$ yield a quasiperiodic orbit and a unique time-stationary measure.   
Next we note:
%To prove that indeed we may restrict to the Gibbs measures when we are interested in the long-time limit of the dynamics we note: 
\begin{enumerate}
\setcounter{enumi}{1}
\item The translation-invariant measures $\mu'$ which have no 
decrease of relative entropy density relative to the invariant measure  $\mu_*'$ (equivalently, to one of the Gibbs measures $\mu'_{\phi}$) in one time step, are necessarily Gibbs measures for the same specification as $\mu'_*$.
\end{enumerate}
In particular translation-invariant and dynamically stationary measures are Gibbs measures for the same potential as $\mu'_*$ and their relative entropy density is also zero. In short: Zero entropic loss implies Gibbsianness. To use such a connection is similar in spirit to the IPS case (where the analogous connection was  termed "Holley's argument" \cite{Ho71,Li85}).  
However, the proof is different
%,  how this is done on a technical level is completely different, 
and in fact we get stronger statements than for the IPS case (compare Theorem \ref{ZeroEntrGibbs} below). We are inspired in this part by the paper \cite{DaLoRo02} about probabilistic cellular automata with strictly local update rules (see also \cite{Ku84}). We also use the relation between, on the one hand,  the decrease of relative entropy density between a general starting measure and the invariant measure  $\mu_*'$ under application of the dynamics, and on the other hand, relative entropies between corresponding time-reversed transition operators. Using K\"unsch's ideas from \cite{Ku84} we then derive the desired DLR equation to identify the Gibbs measures. Technically there are also differences in our treatment to these papers: We need to take proper care of non-localities, but we are able to bypass  a K\"unsch-type representation of the transition operator in terms of double-Gibbs potentials and work directly with specifications, taking advantage of their properties we have to our disposition in our case, which considerably simplifies things.  (For background on how to go from specifications to potentials see \cite{Su73,Ko74,Ku01}.)

The non-reversible time-evolutions we consider here suggest another set of questions, namely whether there are any non-Gibbsian pathologies 
along the trajectories depending on starting measures as found for reversible dynamics in \cite{EnFeHoRe02,EnRu09,KuRe06,KuNy07,ErKu10,EnFeHoRe10,EnKuOpRu10,JaKuRuWe14,FeHoMa13a}.

\bigskip
\textbf{Acknowledgement: }  
This work is supported by the Sonderforschungsbereich SFB $|$ TR12-Symmetries and Universality in Mesoscopic Systems. Christof K\"ulske thanks Universit\'{e} Paris Diderot - Paris 7 for kind hospitality and Giambattista Giacomin for stimulating discussions.

\section{Equilibrium model and rotation dynamics}
In this section we present the equilibrium model also exhibited in \cite{JaKu12}. Further we define the updating via the three-step procedure {\em Sample-Rotate-Project}  and show locality properties.
\subsection{The equilibrium model}\label{The equilibrium model}
\textit{The first-layer model: }We have to first introduce a continuous-spin model which is given in terms of a Gibbsian specification for an absolutely summable Hamiltonian acting on lattice-configurations with continuous local state space. More precisely we consider an $S^1$-rotation invariant and translation-invariant Gibbsian specification $\g^\Phi$ on the lattice $G=\Z^d$, with local state space $S^1=[0,2\pi)$. Let this specification $\g^\Phi=(\g^\Phi_\L)_{\L\subset G}$ 
be given in the standard way by an absolutely summable, $S^1$-invariant and translation-invariant potential 
$\Phi=(\Phi_A)_{A\sb G, A\text{ finite}}$, w.r.t to the Lebesgue measure $\l$ on the spheres. 
This means that the Gibbsian specification is given by the family of probability kernels
\begin{equation*}\label{First_Layer_Specification*}
 \begin{split}
\g^\Phi_\L(B|\eta)=\frac{\int 1_B(\s_\L\eta_{\L^c}) \exp(-H_\L(\s_\L\eta_{\L^c}))\l^{\otimes\L}(d\s_\L)}{\int \exp(-H_\L(\s_\L\eta_{\L^c}))\l^{\otimes\L}(d\s_\L)}
\end{split}
\end{equation*}
for finite $\L\sb G$ and Hamiltonian $H_\L=\sum_{A\cap\L\neq\emptyset}\Phi_A$ applied to a measurable set $B\sb(S^1)^G$ and a boundary condition $\eta\in(S^1)^G$ (for details on Gibbsian specifications see \cite{Ge11}). We use notation $\L^c:=G\setminus\L$. 
%$H_\L$ also has to be differentiable under variation at a single site and these partial derivatives have to be uniformly bounded. 
A standard example of such a model is provided by the nearest-neighbor scalarproduct interaction rotator model with Hamiltonian 
\begin{equation}\label{Metric_Family}
H_{\Lambda}(\s_\L\eta_{\L^c}) = -\beta \sum_{i,j\in \Lambda: i\sim j  }\cos(\s_i-\s_j) -
 \beta \sum_{i\in \Lambda,j \in \Lambda^c: i\sim j  } \cos(\s_i-\eta_j).
\end{equation}

Denote by  $\GG(\g^\Phi)$ the simplex of the Gibbs measures corresponding to this specification,  
which are the probability 
measures $\mu$ on $(S^1)^G$ which satisfy the DLR-equation $\int\mu(d\eta) 
\g^\Phi_\L(B|\eta)=\mu(B)$ for all finite $\L$. Denote by  $\GG_{\theta}(\g^\Phi)$ the lattice
translation-invariant Gibbs measures.

We will make as an assumption on the class of potentials (Hamiltonians)
we discuss moreover that it has a continuous symmetry breaking in the following sense. Assume that the extremal translation-invariant Gibbs measures can be obtained as weak limits with 
homogeneous boundary conditions, i.e with $\eta_\phi\in(S^1)^G$ defined as $(\eta_\phi)_i=\phi$ for all $i\in G$ and $\phi \in S^1$ we have $$\text{ex } \GG_{\theta}(\g^\Phi)=\{ \mu_\phi |\mu_\phi = \lim_{\L\nearrow G} \g^{\Phi}_{\L}(\cdot|\eta_\phi) , \phi \in S^1\}.$$ 
We further assume that
%\begin{equation}\label{rotation_generator}
% \begin{split}
%\text{ex } \GG_{\theta}(\g^\Phi)=\{ \mu_\phi |\mu_\phi = \lim_{\L\nearrow G} \g^{\Phi}_{\L}(\cdot|\eta_\phi) , \phi \in S^1\}% \xrightarrow{\sim}  S^1
%\end{split}
%\end{equation}
different boundary conditions $\eta_\phi$ yield different measures so that there is a  unique labelling of states $ \mu_\phi $ 
by the angles $\phi$ in the sphere $S^1$. It is a non-trivial proven fact that  
this assumption is true in the case of the standard rotator model \eqref{Metric_Family} in $d=3$  for $\l$-a.a temperatures in the low-temperature region as discussed in \cite{FrSiSp76,MaSh11,Pf82}. Here 
the unique labelling can be given by the local magnetization $\mu_\phi(\s_0)=m e_\phi$ where $1>m>0$ is the temperature-dependent length of the unit vector $e_\phi\in S^1$ with angle $\phi$.

\medskip
\textit{The second-layer model: }We will now describe the discretization transformation which maps the continuous-spin model 
to a discrete-spin model.
Denote by $T$ the local coarse-graining with equal arcs, i.e
$T:[0,2\pi)\mapsto \{1, \dots, q\}$ where $T(\phi):=k$ iff $2\pi(k-1)/q\leq\phi<2\pi k/q$.
%of $[0,2\pi)$ to $\{1, \dots, q\}$, i.e. with $S_k:=2\pi/q[k-1,k)$, $S^1=\bigcup_{k=1}^{q}S_k$ and $T(s)=k$ if $s\in S_k$. 
Extend this map to infinite-volume configurations 
by performing it sitewise. We will refer to the image space $\O:=\{0,\dots,q-1\}^G$ as the
coarse-grained layer. 
In particular we will consider images of infinite-volume measures under $T$. 

We will need to choose the parameter of this discretization $q\geq q_0(\Phi)$ large enough so that the image measures of first-layer Gibbs measures are again Gibbs measures for a discrete specification on the coarse-grained layer. 
That this is always possible follows from the earlier works \cite{KuOp08,EnKuOp11}. 
More precisely, we assume that the condition from Theorem 2.1 of 
\cite{EnKuOp11} is fulfilled 
(ensuring a regime where the Dobrushin uniqueness condition holds 
for the so-called constrained first-layer models where the Dobrushin condition is a weak dependence condition implying uniqueness and locality properties). 
Note, as in our notation the usual inverse temperature parameter $\b$ is incorporated into $\Phi$, for $\b$ tending to infinity so does $q_0(\Phi)$.

To talk about the correspondence between the continuous and the discrete system we need 
to make explicit the relevant Gibbsian specification for the latter.
To do so define a family of kernels $\g'=(\g'_\L)_{\L\subset G, \L \text{ finite}}$ for the discretized model by   
\begin{equation}\label{Coarse_Specification}
 \begin{split}
\g'_{\L} (\s'_{\L} | \s' _{\L^c})
&=\frac{\int\mu_{\L^c}[\s'_{\L^c}](d\s_{\L^c})\int\l^{\otimes\L}(d\s_\L)e^{-H_{\L}(\s_\L\s_{\L^c})}
1_{T(\s_\L)=\s'_{\L}}}{\int\mu_{\L^c}[\s'_{\L^c}](d\s_{\L^c})\int\l^{\otimes\L}(d\s_\L)e^{-H_{\L}(\s_\L\s_{\L^c})}
}
\cr 
&=\frac{\mu_{\L^c}[\s'_{\L^c}](\l^{\L} (e^{-H_{\L}}
1_{\s'_{\L}}) ) 
}{\mu_{\L^c}[\s'_{\L^c}](\l^{\L} (e^{-H_{\L}}))}
\cr 
\end{split}
\end{equation}
where the second line is a short notation that we want to adapt in the sequel. The measure $\mu_{\L^c}[\s'_{\L^c}]$ is the unique continuous-spin Gibbs measure for a system on the smaller volume $\L^c$ with conditional specification obtained by deleting all interactions with $\L$ and constrained to take values $\s_{\L^c}$  with discretization 
images $T(\s_{\L^c})=\s'_{\L^c}$.  
For more details and definition of $\mu_{\L^c}[\s'_{\L^c}]$ in terms of formulas 
see \cite{JaKu12} Section 2. 
Note that these constrained Gibbs measures are well-defined and have nice locality properties for sufficiently fine discretization $q\geq q_0(\Phi)$, see \cite{EnKuOp11,KuOp08} and below. For general background on constrained Gibbs measures in the context of preservation of Gibbsianness see \cite{EnFeSo93,Fe05,KuLeRe04}.  We note that $\g'$ is indeed a quasilocal specification and the discretized Gibbs measures are Gibbs for $\g'$. 

\medskip
\textit{The relation between first- and second-layer models: }
The infinite-volume discretization map $T$ is injective when applied to the set of translation-invariant extremal Gibbs states in the continuum model $\text{ex } \GG_{\theta}(\g^\Phi)$. 
More precisely we have the following proposition proved in \cite{JaKu12}.
\begin{prop}\label{Bijection} Let $q\geq q_0(\Phi)$, then $T$ is a bijection from $\text{ex } \GG_{\theta}(\g^\Phi)$ to $\text{ex } \GG_{\theta}(\g')$ 
with inverse given by the kernel % $\g_G^{\s'}(d\s)$. 
$\mu_G[\s'](d\s)$.
\end{prop}  
Here $\mu_G[\s'](d\s)$ is the unique conditional continuous-spin Gibbs measure on the whole volume $G$. 
%It is important to understand that this kernel gets us back from a discrete-spin Gibbs measure to a continuous-spin Gibbs measure in a way which does not depend on the choice of the initial measure. 
%This is crucial for the possibilityto construct a rotation generator $L$ with the desired properties, as we will see.
%The fact that $T \mu :=\mu\circ T^{-1}$ is Gibbs for $\g'$ when $\mu$ is Gibbs for $\g^\Phi$, is already proved in \cite{EnKuOp11,KuOp08} and based on the uniform Dobrushin condition 
%on the coarse-graining. The part that each translation-invariant discrete Gibbs measure has a discretization preimage in the continuous Gibbs measures 
%is new and uses the Gibbs variational principle which involves considerations of 
%relative entropy densities (see \cite{Ge11}).
Regarding part 1 of Theorem \ref{PCA_Theorem} we have the following corollary.
\begin{cor}\label{QuasilocalBijection}
Let $\Phi$ be in the phase-transition region, $q\geq q_0(\Phi)$ and $m=|\mu(\s_0)|$ the uniform local magnetization length for all $\mu\in\text{ex } \GG_{\theta}(\g^\Phi)$. Under these assumptions, the mapping $\psi: \O\mapsto S^1$
\begin{equation*}\label{psi}
\begin{split}
\psi(\s'):=\mu_G[\s'](\s_0)/m
\end{split}
\end{equation*}
is quasilocal and the expectation map 
$\nu \mapsto \nu(\psi)$ is a measurable bijection from $\text{ex } \GG_{\theta}(\g')$ to the sphere $S^1\subset \R^2$. 
\end{cor}
\textbf{Proof: }
The quasilocality of $\psi$ follows from Dobrushin uniqueness arguments given in \cite{JaKu12} for $q\geq q_0(\Phi)$. By the continuous symmetry breaking of the first layer model and the surjectivity of $T$, for every $\phi\in S^1$ there exists a $\mu'_\phi\in\text{ex } \GG_{\theta}(\g')$. Further we have
\begin{equation*}\label{psi_bij}
\begin{split}
\mu'_\phi(\psi)=\frac{1}{m}\int_\O\mu'_\phi(d\s')\mu_G[\s'](\s_0)=\frac{1}{m}\mu_\phi(\s_0)=e_\phi
\end{split}
\end{equation*}
and different $\phi$ yield different $\mu'_\phi$. 
$\Cox$
%For more details on the equilibrium models we refer to \cite{JaKu12}.

\subsection{The updating mechanism}\label{The updating mechanism}

%Take a suitable fixed time $\t\in [0,2 \pi]$. 
%We will play with this time $\t$ in a moment, but let us be general in the beginning. 
We define an updating at finite time $0\leq\t<2\pi$, according to the {\em Sample-Rotate-Project} algorithm, as follows: Given a discrete-spin configuration $\s'$ we perform the following steps: 
\begin{enumerate}
\item Sample a continuous-spin configuration $\s$ according to the conditional measure $\mu_G[\s'](d\s)$.
% where $\s'$ is a coarse-grained configuration.
%conditional on $\s'$, describing an infinite configuration of boxes 
\item Rotate deterministically the resulting continuous-spin configuration $\s \mapsto \s +\t$ jointly in all sites by the same angle $\t$.  
\item Project the rotated configuration using the discretization map $T$, i.e look at the coarse-grained configuration $T(\s+\t)$.  
\end{enumerate}

The resulting kernel on discrete spins, describing the probability distribution of $T(\s+\t)$ where $\s$ is distributed according to $\mu_G[\s'](d\s)$ for a given initial configuration $\s'$, we denote by $M_\t(\s', \cdot\,)$. This \textit{transition operator} can be expressed via
\begin{equation}\label{OneStepUpdate}
\begin{split}
M_\t(\s',\eta'_\L)=\mu_G[\s'](\eta'_{\L,\t})
%=\frac{\mu_{G\ba\L}[\s'_{G\ba\L}](\l^{\L}(e^{-H_\L}1_{\s'_\L}1_{\eta'_{\L,\t}}))}{\mu_{G\ba\L}[\s'_{G\ba\L}](\l^{\L}(e^{-H_\L}1_{\s'_\L}))}
\end{split}
\end{equation}
where $\L$ is a finite set of sites and 
\begin{equation*}\label{Eta_tau}
\begin{split}
\eta'_{\L,\t}:=T^{-1}(\eta'_\L)-\t \in (S^1)^\L. 
\end{split}
\end{equation*}
In words: $\eta'_{\L,\t}$ is obtained by joint rotation by $-\t$ of the segments of the sphere prescribed by $\eta'_\L$.
It will be convenient to use the rewriting
\begin{equation}\label{rewriting}
\begin{split}
M_\t(\s',\eta'_\L)%=\mu_G[\s'](\eta'_{\L,\t})
=\frac{\mu_{\L^c}[\s'_{\L^c}](\l^{\L}(e^{-H_\L}1_{\s'_\L}1_{\eta'_{\L,\t}}))}{\mu_{\L^c}[\s'_{\L^c}](\l^{\L}(e^{-H_\L}1_{\s'_\L}))}
\end{split}
\end{equation}
which follows from reorganization of terms in the Hamiltonian in the Dobrushin uniqueness regime, see \cite{JaKu12} Section 2.

%the coarse-grained configuration $\eta'_\L$ interpreted as a configuration of sphere-segments jointly rotated by the angle $-\t$, i.e $\eta'_{\L,\t}:=T^{-1}(\eta'_\L)-\t=\times_{i\in\L}[\eta_i'2\pi/q-\t,(\eta_i'+1)2\pi/q-\t)\subset[0,2\pi)^\L $.
%$\eta'_{i,\t}:=[\eta_i'|^l-\t,\eta_i'|^r-\t]\subset[0,2\pi)$. 

%\medskip
%Intuitively this representation can best be understood by reading the Sample-Rotate-Project backwards. The probability to see a certain configuration $\eta'_\L$ under $M_\t(\s',\dots)$ 
%%for $0\leq\t<2\pi/q$. In this case the probability to 

%see a configuration which is rotated up one step inside the volume $\D\sb\L$ and 

\medskip
Notice that $M_\t(\s',\eta'_\L)=0$ if and only if $\s'_{\L,0}\cap\eta'_{\L,\t}=\emptyset$. If $M_\t(\s',\eta'_\L)>0$ we call the configuration $\eta'$ \textit{accessible for $\s'$ in $\L$}. For $0\leq\t\leq2\pi/q$ the updating has indeed Bernoulli-increments since in this case the numerator of \eqref{rewriting} is zero for $\eta'_\L\neq\s'_\L+\{0,1\}^\L$.

\medskip
The following two propositions verify the locality properties addressed in Theorem \ref{PCA_Theorem}.
%The transition operator is not a PCA in the strict sense since it is not a product measure. Still it can be understood as a \textit{weak PCA} in the sense described by the following proposition.
%\bigskip
%The Markov transition kernel $M_\t$ works well on the discrete spin Gibbs measures $\mu'_{\phi}$: 
\begin{prop}\label{Mixing_Past}  
For any $\a_2>0$ we can choose the descretization $q=q(\a_2)$ fine enough such that the update kernel $M_\t(\s', d\eta')$ is uniformly spatially mixing in the future with exponent $\a_2>0$. 
\end{prop} 

{\bf Proof: } Let $\s'$ be any starting configuration. Further let $A$ be a set of configurations measurable w.r.t the sigma-algebra corresponding to a finite volume $\L$ and $B$ be a set of configurations measurable w.r.t the sigma-algebra corresponding to a volume $\D^c$ where $\D$ is a finite volume. Further assume 
%$B$ to contain at least one accessible configuration $\eta'$ for $\s'$ in $\bar\L\subset\D^c$ and hence 
$M_\t(\s',B)>0$, then we have
% finite volumes with $\L\subset\D$ and $\bar\L\subset\D^c$ we have
%\begin{equation*}\label{in the future2}
% \begin{split}
%|M_\t(\o',\eta'_\L|\eta'_{\bar\L}) - M_\t(\o',\eta'_\L)| \leq K_2 \sum_{i\in \L, j\in \D^c }e^{- \a_2 |i-j| }
%\end{split}
%\end{equation*}
%\begin{equation*}\label{in the future2}
% \begin{split}
%|M_\t(\o',A|B) - M_\t(\o',A)| \leq K_2 \sum_{i\in \L, j\in \D^c }e^{- \a_2 |i-j| }
%\end{split}
%\end{equation*}
\begin{equation*}\label{zwei}
 \begin{split}
M_\t(\s',A|B) - M_\t(\s',A)
&=\mu_G[\s'](B_{\t})^{-1}\mu_G[\s'](A_\t\cap B_\t)-\mu_G[\s'](A_{\t}).\cr
\end{split}
\end{equation*}
where we used notation analogue to \eqref{OneStepUpdate} for $A_\t, B_\t$.
Notice, 
%with $K_2$ some positive constant. Notice $M_\t(\o',\eta'_{\bar\L})>0$ and using notation \eqref{OneStepUpdate}
%such that  at the following decay estimate holds: For all ,  for all sets $A$ in the sigma-algebra measurable w.r.t. the sigma-algebra generated by the $\eta$-coordinates in $\L$, and sets $B$ in the sigma-algebra generated by the $\eta$-coordinates in $G\setminus\D$ one has 
%\begin{equation*}\label{zwei}
% \begin{split}
%M_\t(\o',\eta'_\L|\eta'_{\bar\L}) - M_\t(\o',\eta'_\L)
%&=\mu_G[\o'](\eta'_{\bar\L,\t})^{-1}\mu_G[\o'](\eta'_{\L\cup\bar\L,\t})-\mu_G[\o'](\eta'_{\L,\t}).\cr
%\end{split}
%\end{equation*}
$\mu_G[\s']$ is uniquely specified by the constrained specification $\g^{\s'}$ (for details see \cite{JaKu12}) and thus we can write
%\begin{equation*}\label{zwei}
% \begin{split}
%\mu_G[\o'](\eta'_{\L\cup\bar\L,\t})=\int_{\eta'_{\bar\L,\t}}\g^{\o'}_{\D}(\eta'_{\L,\t}|\o)\mu_G[\o'](d\o).\cr
%\end{split}
%\end{equation*}
\begin{equation*}\label{zwei}
 \begin{split}
\mu_G[\s'](A_\t\cap B_\t)=\int_{B_{\t}}\g^{\s'}_{\D}(A_{\t}|\s)\mu_G[\s'](d\s).\cr
\end{split}
\end{equation*}
$\g^{\s'}$ is in the Dobrushin uniqueness region with Dobrushin matrix $\bar C$ uniformly in $\s'$ thus
%, i.e  $1>\bar c:=\sup_i \sum_j \bar C_{ij}$ 
%(for more details see \cite{JaKu12} Section 2). 
by Theorem 8.23 (ii) in \cite{Ge11} we have
\begin{equation*}\label{zwei}
 \begin{split}
%\sup_{\eta'_{\bar\L,\t}}\Vert\g^{\o'}_{\D}(\eta'_{\L,\t}|\cdot)-\mu_G[\o'](\eta'_{\L,\t})\Vert\leq\sum_{i\in \L, j\in \D^c }\bar D_{ij}
\sup_{A_\t}\Vert\g^{\s'}_{\D}(A_\t|\cdot)-\mu_G[\s'](A_\t)\Vert\leq\sum_{i\in \L, j\in \D^c }\bar D_{ij}
\end{split}
\end{equation*}
%where the supremum runs over all sets $A$ measurable w.r.t the sigma-algebra corresponding to the volume $\L$. 
with $\bar D:=\sum_{n\geq0} \bar C^n$. Now we can pick $q\geq q(\a_2)$ large enough such that the Dobrushin matrix has exponential decay (see \cite{Ge11} Remark 8.26 and \cite{JaKu12} Lemma 2.8) and hence $\bar D_{ij}\leq K_2 e^{- \a_2 |i-j| }$ which finishes the proof.
$\Cox$

\begin{prop}\label{Mixing_Futur}  
The update kernel $M_\t$ satisfies the property of exponential locality from the past with exponent $\a_1$. 
\end{prop} 

{\bf Proof: } Using notation as in \eqref{OneStepUpdate} we have 
\begin{equation*}\label{Product}
 \begin{split}
\d_{i}( M(\cdot ,\eta'_\L))=\sup_{{\s',\tilde \s':}\atop{\s'_{i^c}=\tilde\s'_{i^c}}}|\mu_G[\s'](\eta'_{\L,\t})-\mu_G[\tilde\s'](\eta'_{\L,\t})|
\leq \sum_{j\in\L}\bar D_{ji}
%K_1 e^{-\a_1 d(i,\L)}.
\end{split}
\end{equation*}
since $\mu_G[\s']$ is the unique Gibbs measure for the specification $\g^{\s'}$ which is in the Dobrushin uniqueness region see Theorem 8.20 in \cite{Ge11}.
%We have to show that there is a positive constant $K_1$ such that 
%%the variation w.r.t. $\o_i$ of the probability $M(\o,\eta_V)$ is bounded by the exponential estimate 
%\begin{equation*}\label{Product}
% \begin{split}
%%&\d_{i}( M(\cdot ,\eta'_\L))= 
%\sup_{{\s',\tilde \s':}\atop{\s'_{i^c}=\tilde\s'_{i^c}}}|\mu_G[\s'](\eta'_{\L,\t})-\mu_G[\tilde\s'](\eta'_{\L,\t})|
%\leq |\L|K_1 e^{-\a_1 d(i,\L)}.
%\end{split}
%\end{equation*}
%%Here we use notation $i^c=G\setminus \{i\}$ and $d(i,V)=\min_{j\in V}d(i,j)$ for some semimetric $d$ on $G$.
As above we can choose $q\geq q(2\a_1)$ to be large enough such that %the Dobrushin matrix is exponentially decaying we have using \cite{Ge11} Remark 8.26 and \cite{JaKu12} Lemma 2.8 $\sum_{j\in\L}
$\bar D_{ji}\leq C_1e^{-2\a_1 |i-j|}$ and thus for $i\notin\L$ and $\a_1$ sufficiently large
\begin{equation*}\label{Product}
 \begin{split}
\sum_{j\in\L}\bar D_{ji}&\leq C_1 \sum_{j\in\L} e^{-2\a_1 |i-j|}
\leq C_1 \sum_{k=|i-\L|}^\infty  e^{-2\a_1 k}\#\{j\in \Z^d, |j-i|=k \}\cr
& \leq C_2 e^{-\a_1 |i-\L|}.
\end{split}
\end{equation*}
%If $i\in\L$: $\sum_{j\in\L} e^{-\a_1 [|i-j|-|i-\L|]}=\sum_{j\in\L} e^{-\a_1 |i-j|}=\sum_{k\in\L-i} e^{-\a_1 |k|}\leq\sum_{k\in G} e^{-\a_1 |k|}\leq C_2 $ and if 
%: $\sum_{j\in\L} e^{-\a_1 [|i-j|-|i-\L|]}\leq\sum_{j\in\L} e^{-\a_1 |j-\L^c|}\leq C_3 $. 
$\Cox$

\bigskip
As already addressed in Section \ref{Comp}, in \cite{JaKu12} we present a rotation dynamics in continuous-time as an IPS. Specifically we exhibit a Markov generator $L$ with the property that for the associated semigroup $(S_t)_{t\geq0}$ we have $S_t\mu'_\phi=\mu'_{\phi+t}$ where $\mu'_\phi\in\text{ex } \GG_{\theta}(\g')$ (see \cite{JaKu12} Theorem 1.3). The exponential  locality property \eqref{Exp_Loc} is a  consequence of Lemma 3.4 in \cite{JaKu12}. 
%Let us explain why the result of JaKu give as the arbitrarily well localized 
%transition rates as promised in Theorem 1.2. 
%Recall the definition of the rates  from JaKu. 
More precisely, given any arbitrarily large decay exponent $\a_1$ 
we may ensure the decay \eqref{Exp_Loc} by choosing $\b$ and $q$ in \cite{JaKu12} 
such that the corresponding Dobrushin-constant $\bar c=\bar c (\b,q)$ (of formula (28) of \cite{JaKu12}) 
is sufficiently small. 
Now choose $\b$ sufficiently large such that the initial rotator models has a phase transition. 
Since $\bar c (\b,q)\downarrow 0 $ for $\b$ fixed, with $q\uparrow \infty$, the statement 
now follows for $q$ large. 

\medskip
In the present discrete-time setting the mapping $\t \mapsto M_\t$ can not be expected to be a semigroup. In particular at finite $\t$ the transition operator $M_\t$ differs from the continuous rotation semigroup $(S_t)_{t\geq0}$. They become equal at infinitesimal $\t$ and this is the route for the definition of $(S_t)_{t\geq0}$ (see the proof of Theorem 1.3 in \cite{JaKu12} Section 3.2). Nevertheless $M_\t$ also possesses the rotation property presented in the following proposition.

\begin{prop}\label{RotationProp} 
Discrete Gibbs measures transform in a covariant way, i.e \linebreak
$M_\t \mu'_{\phi}= \mu'_{\phi+\t}$ for all $\mu'_\phi\in\text{ex } \GG_{\theta}(\g')$ and $0\leq\t<2\pi$.
\end{prop}
In particular $\II:=\GG_{\theta}(\g')$ is the set of translation-invariant measures where the discrete-time process acts like a rotation. Applying the function $\psi$ from Section \ref{The equilibrium model} this proves part 1 and 2 in Theorem \ref{PCA_Theorem}.

\bigskip
\textbf{Proof:} It suffices to prove that for all $\mu'_\phi\in\text{ex } \GG_{\theta}(\g')$ and discrete-spin test-functions $f$ the following equality holds: 
\begin{equation*}\label{zwei}
 \begin{split}
 &\int \mu'_\phi(d\s')\int M_\t(\s',d\eta')f(\eta')=\int\mu'_{\phi+\t}(d\eta')f(\eta')
\end{split}
\end{equation*}
By Proposition \ref{Bijection} the conditional probability under coarse-graining of a continuous-spin measure ordering in direction $\phi$ does not depend on $\phi$, i.e for continuous-spin test-functions $g$ we have 
\begin{equation}\label{Kernel}
\begin{split}
&\int \mu_\phi(d\s)g(\s)=\int\mu'_{\phi}(d\s')\mu[\s'](d\s)g(\s).
\end{split}
\end{equation}
%The definition of $\mu'_{\phi}$ says 
%\begin{equation}\label{vier}
% \begin{split}
%	 &\int\mu'_\phi(d\s')f(\s')=\int\mu_{\phi}(d\s)f(T\s).
%\end{split}
%\end{equation}
Thus we can write 
\begin{equation*}\label{vier2}
 \begin{split}
\int \mu'_\phi(d\s')\int M_\t(\s',d\eta')f(\eta')&=\int \mu'_\phi(d\s')\int \mu[\s'](d\s)f(T(\s+\t) )\cr
 &=\int \mu_\phi(d\s)f(T(\s+\t) )\cr
 &=\int \mu_{\phi+\t}(d\s)f(T(\s) )\cr
 &=\int\mu'_{\phi+\t}(d\eta')f(\eta')\cr
\end{split}
\end{equation*}
where the first equality is the definition of $M_\t$, the second equality is the independence of conditional probability on $\phi$ (see \eqref{Kernel}), the third equality is the property of the continuous-spin Gibbs measures to transform under rotations and the last inequality is the definition of the coarse-grained measures. 
$\Cox$

\begin{figure}[h]
$$
\begin{xy}
  \xymatrix{
      \PP_\theta((S^1)^{\Z^3})\supset \ar[dd]^T & \text{ex } \GG_{\theta}(\g^\Phi) \ar[rrrr]^{\mu\mapsto R_\t\mu} \ar@/_0,5cm/[dd] & & & & \text{ex } \GG_{\theta}(\g^\Phi) \ar@/_0,5cm/[dd]  \\ \\
      \PP_\theta((\Z_q)^{\Z^3})\supset \ar[dd]^{\mu'\mapsto \mu'(\psi)} & \text{ex } \GG_{\theta}(\g') \ar[rrrr]^{\mu'\mapsto \int\mu'(d\s')M_\t(\s',\cdot)} \ar@/_0,5cm/[dd] \ar@/_0,5cm/[uu] & & & &   \text{ex } \GG_{\theta}(\g') \ar@/_0,5cm/[dd] \ar@/_0,5cm/[uu]\\ \\
      \R^2\supset & S^1 \ar[rrrr]^{\phi\mapsto\phi+\t} \ar@/_0,5cm/[uu] & & & &   S^1 \ar@/_0,5cm/[uu]  
  }
\end{xy}
$$
\caption{\scriptsize{The discretization map $T$ and labelling map $\mu'\mapsto \mu'(\psi)$ become bijections when applied to the extremal translation-invariant Gibbs measures. The transition kernel $M_\t$ reproduces the deterministic rotation actions $R_\t$ and $\phi\mapsto\phi+\t$.}}
\label{Diagram}
\end{figure}
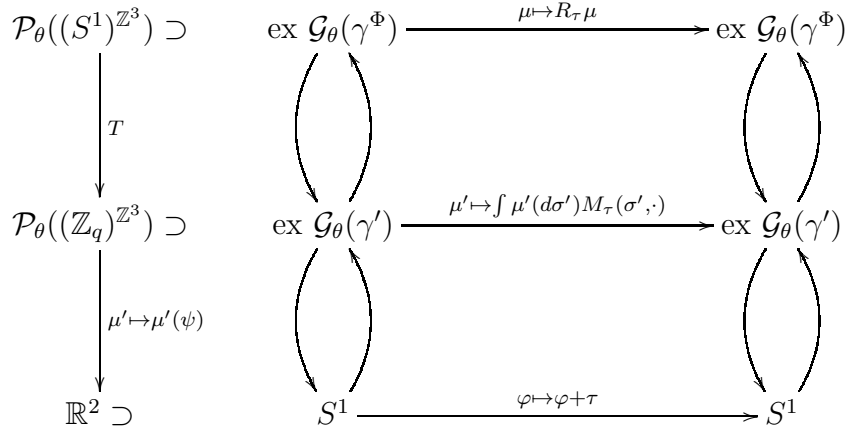

\bigskip
To summarize the interplay between the discretization and the dynamics let us consider the joint rotation of a first-layer measure $\mu$ by an angle $\t$ written as $R_\t\mu$. Then, under $R_\t$ and the Markov transition operator $M_\t(\s',d \eta')$, the diagram in Figure \ref{Diagram} is commutative.
%\begin{figure}[h]
%$$
%\begin{xy}
%  \xymatrix{
%       (S^1)^{\Z^d} \ar[rrr]^{\o\mapsto R_\t\o} \ar@/_0,5cm/[dd]_T &  &  &  (S^1)^{\Z^d} \\ \\
%      \{1,\dots,q\}^{\Z^d} \ar[rrr]_{\o'\mapsto M_\t(\o',d \hat \o')
%            } \ar@/_0,5cm/[uu]_{\o\sim \mu_G[\o'](d\o)}        &   & &  \{1,\dots,q\}^{\Z^d}
%  }
%\end{xy}
%$$
%\caption{\scriptsize{Definition of transition kernel $M_\t$ for PCA
%}
%}
%\label{Diagram}
%\end{figure}
%\begin{figure}[h]
%$$
%\begin{xy}
%  \xymatrix{
%      \text{ex } \GG_{\theta}(\g^\Phi) \ar[rrr]^{\mu\mapsto R_\t\mu} \ar@/_0,5cm/[dd]_T &  &  & \text{ex } \GG_{\theta}(\g^\Phi) \ar[dd]^T  \\ \\
%      \text{ex } \GG_{\theta}(\g') \ar[rrr]_{\mu'\mapsto \int\mu'(d\s')M_\t(\s',\cdot)} \ar@/_0,5cm/[uu]_{\mu'\mapsto\int\mu'(d\o')\mu_G[\o'](d\o)}        &   & &   \text{ex } \GG_{\theta}(\g')   
%  }
%\end{xy}
%$$
%\caption{\scriptsize{Equivariance property of the bijective discretization map $T$ for the deterministic rotation action $R_\t$ and the action of the transition kernel $M_\t$.}}
%\label{Diagram}
%\end{figure}

\medskip
In order to further explain the transition kernel $M_\t$, let us consider the case of the zero initial continuous-spin Hamiltonian
%, i.e. $\Phi_A=0$ for all $A\subset G$ 
as an example. In this case the spins become independent and it is sufficient to consider a single spin with the Gibbs distribution 
being the uniform distribution on the circle $\s_0\sim \hbox{Unif}([0,2\pi])$. The image measure under $T$ of the continuous-spin Gibbs measure is again the uniform distribution on the discretized circle $\s'_0\sim \hbox{Unif}(\{1,\dots,q\})$. 
%
%We have, conditional on the discretized spin $\o'_0\in \{1,\dots, q\}$, 
%$\s^{\o'_0}\sim \hbox{Unif}([\o'-1,\o']\frac{2\pi}{q})$. 
Hence, for $0\leq\t<2\pi/q$ we have that 
\begin{equation*}\label{zwei}
 \begin{split}
M_\t(\s'_0,\{\s'_0,\s'_0+1\})=1\hspace{0.2cm}\text{ and }\hspace{0.2cm}M_\t(\s'_0,\s'_0+1)=\t\frac{q}{2\pi}.\cr
 \end{split}
\end{equation*}
This describes a random walk $(\s'_0(n\t))_{n\in\N}$ on the discretized circle $\{1,\dots,q\}$ which can only move up by one step with probability proportional to $\t$ or stay where it is. Clearly for larger parameters $\t$ it is more probable to jump. In that sense $\t$ controls the velocity of the walker.  

Taking the limit $\t\downarrow 0$ the process converges to a Poisson process which moves around the discretized circle where jumps are made with intensity $q/2\pi$.

\bigskip
Finally let us point out that the uniformly mixed Gibbs measure $\mu_*':=\frac{1}{2\pi}\int_0^{2\pi} d\phi \mu'_{\phi}$ with $\mu'_\phi\in\text{ex } \GG_{\theta}(\g')$ is time-stationary for $M_\t$. Indeed we have
\begin{equation}\label{MixedMeasure}
 \begin{split}
M_\t\mu_*'=\frac{1}{2\pi}\int_0^{2\pi} d\phi M_\t\mu'_{\phi}=\frac{1}{2\pi}\int_0^{2\pi} d\phi\mu'_{\phi+\t}=\frac{1}{2\pi}\int_0^{2\pi} d\phi\mu'_{\phi}=\mu_*'.
  \end{split}
\end{equation}
In other words, there exists a translation-invariant Gibbs measure that is also time-stationary for the dynamics.
In the following section we prove that any translation-invariant and time-stationary measure must be a Gibbs measure for the same specification.

%{\bf Time-invariance of uniform mixture $\mu^*$ for any $\t$:} From $M_\t \mu'_\phi =\mu'_{\phi+\t}$ we have that the uniform mixture
%$\mu^*:=\frac{1}{2\pi}\int d\phi \mu'_{\phi}$ is invariant for any $\t$. If $\frac{\t}{2\pi}$ is rational there are periodic orbits and we have uncountably many invariant measures sitting on the possible orbits. If $\frac{\t}{2\pi}$ is irrational, the trajectory $k\mapsto M^k_{\t}\mu_{\phi}=\mu_{\phi+ k\t}$ 
%is aperiodic and dense.  

\section{Entropic loss and reversed transition operators}
%\subsection{Properties of the transition kernel $M_\t$ }
%Note that $M_\t$ (and $M$) is not in a strict-sense a PCA, that is 
%the probability $M(\s',\eta'_V)$ to see a finite volume configuration $\eta'_V$ given the 
%initial configuration $\s'$ does {\em not} factorize in the form $M(\s',\eta'_V)=\prod_{i\in V}M(\s',\eta'_i)$.
%
%$M$ is however {\em jointly Gibbsian} in the sense of K\"unsch, SPA1984. 
%This is a consequence of our construction which uses fine enough discretization 
%such that Dobrushin uniqueness holds for the constrained measures.  
%
%Note that $M_\t$ wouldn't be jointly Gibbs because of lack of non-nullness. 
%The measure on the $n's$ would be though, maybe this is enough for 
%the use of K\"unsch which is to follow. CHECK! 
%
We want to connect the two properties, having zero entropy loss under time evolution and being Gibbs w.r.t the same specification as $\mu'_*$. In order to do this relative entropy densities will be considered. More precisely we want to employ an extension to our Gibbsian updating mechanism of arguments which are carried out in \cite{DaLoRo02} Proposition 2.1 and Proposition 2.2 involving time-reversed transition operators. From now on we will always consider translation-invariant measures. 

Let us introduce some notation. For infinite-volume probability measures $\nu',\mu'\in\PP(\{0,\dots,q-1\}^G)$ and a finite set of sites $\L$ the \textit{local relative entropy} is defined as 
\begin{equation*}\label{LRE}
\begin{split}
h_\L(\nu'|\mu'):=\sum_{\s'_\L\in\{1,\dots,q\}^\L}\nu'(1_{\s'_\L})\log\frac{\nu'(1_{\s'_\L})}{\mu'(1_{\s'_\L})}.
\end{split}
\end{equation*}
If $\mu'\in\GG_{\theta}(\g')$, existence of the \textit{specific relative entropy}
\begin{equation*}\label{SRE}
\begin{split}
h(\nu'|\mu'):=\limsup_{\L\uparrow G}\frac{1}{|\L|}h_\L(\nu'|\mu')
\end{split}
\end{equation*}
is guaranteed, where $\L$ varies over hypercubes centered at the origin. Similarly one can define the \textit{specific relative entropy between transition operators w.r.t a base measure $\nu'$} by 
\begin{equation*}\label{SRETrans}
 \begin{split}
 &\HH_{\nu'}(M|\tilde M)=\int \nu'(d\s')\limsup_{\L}\frac{1}{|\L|}h_\L\big( M(\s',\cdot)|\tilde M(\s',\cdot)\big).
   \end{split}
\end{equation*}
Let us define the \textit{joint two-step distribution} $Q_{\nu'}(d\s',d\eta'):=M_\t(\s',d\eta')\nu'(d\s')$. We will consider different conditionings of $Q_{\nu'}$. To keep notation reasonably simple we set the convention and write \textit{$\s'$ for the present configuration} and \textit{$\eta'$ for the future configuration}, just as in $Q_{\nu'}(d\s',d\eta')$. With this we have $Q_{\nu'}(\eta'_\L|\s')=M_\t(\s',\eta'_\L)$ which is independent of $\nu'$. The \textit{backwards transition operator} is given by 
\begin{equation*}\label{BackwardsTrans}
 \begin{split}
\hat M_{\t,\nu'}(\eta',\s'_\L):=Q_{\nu'}(\s'_\L|\eta'). 
   \end{split}
\end{equation*}
Using the short notation $M_\t\nu'(\eta'_\L)=\int\nu'(d\s')M_\t(\s',\eta'_\L)$, the backwards transition operator is characterized by the requirement that 
\begin{equation}\label{BackwardsTransChar}
 \begin{split}
\int\nu'(d\s')\int M_\t(\s',d\eta')f(\s',\eta')=\int M_\t\nu'(d\eta')\int\hat M_{\t,\nu'}(\eta',d\s')f(\s',\eta') 
   \end{split}
\end{equation}
%if \eqref{BackwardsTransChar} 
holds for all local test-functions $f$.  
\begin{lem}\label{GibbsBackwards}
The backwards transition operator for any translation-invariant Gibbs measure 
%$\mu'\in\GG_{\theta}(\g')$ 
is given by $M_{-\t}$ where $M_{-\t}$ is obtained from formula \eqref{OneStepUpdate} for negative $\t$.
\end{lem}
\textbf{Proof: }Let us first check \eqref{BackwardsTransChar} for the extremal Gibbs measure, i.e let $\mu'_\phi\in\text{ex }\GG_{\theta}(\g')$, then 
%
%It is simple to verify that the time reversed transition operator for $M_{\t}$, defined for a starting measure
% $\nu$ to be 
%\begin{equation}\label{zwei}
% \begin{split}
% &\int \nu(d\eta')\int \mu[\eta'](d\eta)f(T(\eta+\t),\eta')=\int (M_\t\nu)(d\s')\int \hat M(d\eta'|\s') f(\s',\eta')
%\end{split}
%\end{equation}
%is given by $\hat M(d\eta'|\s')=M_{-\t}(d \eta' | \s')$ for all starting measures 
%$\mu'_\phi\in\text{ex } \GG_{\theta}(\g')$ is $M_{-\t}$, 
independently of the angle $\phi$ we have
%. This is clear since
%The proof is clear since 
\begin{equation*}\label{BackwardsGibbs}
 \begin{split}
\int \mu'_\phi(d\s')\int M_\t(\s',d\eta')f(\s',\eta')&=\int \mu'_\phi(d\s')\int \mu_G[\s'](d\eta)f(\s',T(\eta+\t))\cr
&=\int \mu'_\phi(d\s')\int \mu_G[\s'](d\eta)f(T(\eta),T(\eta+\t))\cr
&=\int \mu_\phi(d\eta)f(T(\eta),T(\eta+\t))\cr
&=\int \mu_{\phi+\t}(d\eta)f(T(\eta-\t),T(\eta))\cr
&=\int \mu'_{\phi+\t}(d\eta')\int \mu_G[\eta'](d\s)f(T(\s-\t),\eta')\cr
&=\int \mu'_{\phi+\t}(d\eta')\int M_{-\t}(\eta',d\s') f(\s',\eta')
\end{split}
\end{equation*}
%Further notice by \eqref{BackwardsGibbs} for the backwards transition operator for an extremal translation-invariant Gibbs measure we have 
%\begin{equation*}\label{BackwardsGibbs2}
% \begin{split}
%\hat Q_{\mu}(\s_{\L}|\eta)=M_{-\t}(\eta,\s_{\L}).
%   \end{split}
%\end{equation*}
%where $\hat M$ is a $q$-fold composition of $M_{-\t}$. 
where we used equation \eqref{Kernel} two times. By linearity of the integrals the above equation also holds for any convex combination of the extremal Gibbs measures and hence for all $\mu'\in\GG_{\theta}(\g')$.
$\Cox$

\bigskip
Next we consider the entropy loss under $M_\t$ which can be expressed in terms of the backward transition operators.
\begin{lem}\label{Extension of Proposition 2.1. from DaPrLoRo02} 
Suppose that $\mu'$ is a translation-invariant Gibbs measure w.r.t the specification $\g'$ and also time-stationary w.r.t $M_\t$. Then, for any translation-invariant measure $\nu'$ the entropic loss 
can be expressed via 
\begin{equation}\label{Lemma1}
 \begin{split}
 &h(\nu'|\mu')-h(M_\t\nu'|\mu' )=\HH_{M_\t\nu'}(\hat M_{\t,\nu'}| M_{-\t}). 
   \end{split}
\end{equation}
%where $\hat M_{\t,\nu'}$, $\hat M_{\t,\mu'}$ are the time-reversals of the transition operator w.r.t to reference measures $\nu'$ resp. $\mu'$ 
%and 
%\begin{equation*}\label{zwei}
% \begin{split}
% &\HH_{\nu'}(M_1|M_2)=\int \nu'(d\eta') \limsup_{\L}\frac{1}{|\L|}S( M_1(\eta', \cdot |_{\L})| M_2(\eta', \cdot |_{\L})  )
%   \end{split}
%\end{equation*}
%denotes the relative entropy density between two transition operators $M_1, M_2$ averaged over the starting configuration w.r.t the reference measure $\nu'$. 
\end{lem}

In \cite{DaLoRo02} this result is proved for the case of a PCA with sitewise independent local updating. Here we extend it to our case of a weak PCA with quasilocal updating. 

\bigskip
{\bf Proof: }%It might be convenient to use the following notation:
%\begin{equation*}\label{OneStepUpdate}
%\begin{split}
%M_\t(\s',\eta'_\L)=\mu[\s'](\eta'_{\L,\t})=\frac{\mu[\s'_{\L^c}](\l^{\L}(e^{-H_\L}1_{\s'_\L}1_{\eta'_{\L,\t}}))}{\mu[\s'_{\L^c}](\l^{\L}(e^{-H_\L}1_{\s'_\L}))}
%\end{split}
%\end{equation*}
%where $\eta_{i,\t}:=[\eta_i|^l-\t,\eta_i|^r-\t]$. Of course $M_\t(\s',\eta'_\L)=0$ if $\s'_\L\cap\eta'_{\L,\t}=\emptyset$.
%Since the proof in \cite{DaLoRo02} uses the locality of the updating, we would first like to derive some bounds on the error we make by replacing our quasilocal transition operator by a local one. 
Let us suppress all primes in the notation.
% whenever it is unambiguous. 
First notice, the error we make by replacing the starting configuration outside some finite volume is of boundary order. Indeed, let $\xi$, whenever it appears, be an arbitrary but fixed configuration in $\{0,\dots,q-1\}^{G}$, then
% for $M_\t(\s,\eta_\L)\neq0$
%For all bounded and $\mathcal F_\L$-measurable functions $f$ we have
\begin{equation}\label{OneStepError}
\begin{split}
\frac{M_\t(\s_{\L}\xi_{\L^c},\eta_\L)}{M_\t(\s,\eta_\L)}=\frac{\frac{\mu_{G\ba\L}[\xi_{G\ba\L}](\l^{\L}(e^{-H_\L}1_{\s_\L}1_{\eta_{\L,\t}}))}{\mu_{G\ba\L}[\xi_{G\ba\L}](\l^{\L}(e^{-H_\L}1_{\s_\L}))}}{\frac{\mu_{G\ba\L}[\s_{G\ba\L}](\l^{\L}(e^{-H_\L}1_{\s_\L}1_{\eta_{\L,\t}}))}{\mu_{G\ba\L}[\s_{G\ba\L}](\l^{\L}(e^{-H_\L}1_{\s_\L}))}}\leq e^{4\sum_{A\cap\L\neq\emptyset,A\cap\L^c\neq\emptyset}\Vert\Phi_A\Vert}=e^{o(|\L|)}
\end{split}
\end{equation}
%\begin{equation}\label{OneStepError}
%\begin{split}
%\frac{M_\t(\s'_{\L}\xi'_{\L^c},f)}{M_\t(\s',f)}=\frac{\frac{\mu[\xi'_{\L^c}](\l^{\L}(e^{-H_\L}1_{\s'_\L}f))}{\mu[\xi'_{\L^c}](\l^{\L}(e^{-H_\L}1_{\s'_\L}))}}{\frac{\mu[\s'_{\L^c}](\l^{\L}(e^{-H_\L}1_{\s'_\L}f))}{\mu[\s'_{\L^c}](\l^{\L}(e^{-H_\L}1_{\s'_\L}))}}\leq e^{4\sum_{A\cap\L\neq\emptyset,A\cap\L^c\neq\emptyset}\Vert\Phi_A\Vert}=e^{o(|\L|)}
%\end{split}
%\end{equation}
where of course we assumed $M_\t(\s,\eta_\L)\neq0$. 

%In order to get from $M_\t$ to $M$ one can iterate the following argument for $M_{2\t}$.
%\begin{equation}\label{TwoStepError}
%\begin{split}
%\frac{M_{2\t}(\s'_{\L}\xi'_{\L^c},\eta'_\L)}{M_{2\t}(\s',\eta'_\L)}&=\frac{\int\mu[\s'_\L\xi'_{\L^c}](d\zeta'_\t)\mu[\zeta'](\eta_{\L,\t}')}{\int\mu[\s'](d\zeta'_\t)\mu[\zeta'](\eta_{\L,\t}')}\cr
%&=\frac{\sum_{a'_\L}\int\mu[\s'_\L\xi'_{\L^c}](d\zeta'_\t)\mu[\zeta'](\eta_{\L,\t}')1_{a'_\L}(\zeta')}{\sum_{a'_\L}\int\mu[\s'](d\zeta'_\t)\mu[\zeta'](\eta_{\L,\t}')1_{a'_\L}(\zeta')}\cr
%&=\frac{\sum_{a'_\L}\int\mu[\s'_\L\xi'_{\L^c}](d\zeta'_\t)\frac{\mu[\zeta'](\eta_{\L,\t}')}{\mu[a'_\L\xi'_{\L^c}](\eta_{\L,\t}')}\mu[a'_\L\xi'_{\L^c}](\eta_{\L,\t}')1_{a'_\L}(\zeta')}{\sum_{a'_\L}\int\mu[\s'](d\zeta'_\t)\frac{\mu[\zeta'](\eta_{\L,\t}')}{\mu[a'_\L\xi'_{\L^c}](\eta_{\L,\t}')}\mu[a'_\L\xi'_{\L^c}](\eta_{\L,\t}')1_{a'_\L}(\zeta')}\cr
%&\leq e^{o(|\L|)}\frac{\sum_{a'_\L}\mu[\s'_\L\xi'_{\L^c}](a'_{\L,\t})\mu[a_\L'\xi'_{\L^c}](\eta_{\L,\t}')}{\sum_{a'_\L}\mu[\s'](a'_{\L,\t})\mu[a_\L'\xi'_{\L^c}](\eta_{\L,\t}')}\cr
%&=e^{o(|\L|)}\frac{\sum_{a'_\L}\frac{\mu[\s'_\L\xi'_{\L^c}](a'_{\L,\t})}{\mu[\s'](a'_{\L,\t})}\mu[\s'](a'_{\L,\t})\mu[a_\L'\xi'_{\L^c}](\eta_{\L,\t}')}{\sum_{a'_\L}\mu[\s'](a'_{\L,\t})\mu[a_\L'\xi'_{\L^c}](\eta_{\L,\t}')}\cr
%&\leq e^{o(|\L|)}\frac{\sum_{a'_\L}\mu[\s'](a'_{\L,\t})\mu[a_\L'\xi'_{\L^c}](\eta_{\L,\t}')}{\sum_{a'_\L}\mu[\s'](a'_{\L,\t})\mu[a_\L'\xi'_{\L^c}](\eta_{\L,\t}')}= e^{o(|\L|)}\cr
%\end{split}
%\end{equation}
Then since $\mu$ is assumed to be time-stationary
%, i.e $M_\t\mu=\mu$ we have 
the local entropic loss can be expressed in terms of the backwards transition operator and some error term of boundary order. Indeed,
\begin{equation*}\label{EntropyProductionAbove}
\begin{split}
&h_{\L}(\nu|\mu)-h_\L(M_\t\nu|\mu)\cr
%&=\sum_{\s_{\L}\in\{1,\dots,q\}^\L}\nu(\s_{\L})\log\frac{\nu(\s_{\L})}{\mu(\s_{\L})}-\sum_{\eta_\L\in\{1,\dots,q\}^\L}M_\t\nu(\eta_\L)\log\frac{M_\t\nu(\eta_\L)}{\mu(\eta_\L)}\cr
&=\sum_{\eta_\L\in\{1,\dots,q\}^\L}M_\t\nu(\eta_\L)\log\frac{\mu(\eta_\L)}{M_\t\nu(\eta_\L)}+\sum_{\s_{\L}\in\{1,\dots,q\}^\L}\nu(\s_{\L})\log\frac{ \nu(\s_{\L})}{\mu(\s_{\L})}\cr
%&=\sum_{\eta_\L}M_\t\nu(\eta_\L)\sum_{\s_{\L}} Q_\nu(\s_{\L}|\eta_\L)\log\frac{ \nu(\s_{\L})}{M_\t\nu(\eta_\L)}\frac{M_\t\mu(\eta_\L)}{\mu(\s_{\L})}\cr
&=\sum_{\eta_\L}M_\t\nu(\eta_\L)\sum_{\s_{\L}}Q_\nu(\s_{\L}|\eta_\L)\log\frac{M_\t(\s_{\L}\xi_{\L^c},\eta_\L)\nu(\s_{\L})}{M_\t\nu(\eta_\L)}\frac{M_\t\mu(\eta_\L)}{M_\t(\s_{\L}\xi_{\L^c},\eta_\L)\mu(\s_{\L})}\cr
&=\sum_{\eta_\L}M_\t\nu(\eta_\L)\sum_{\s_{\L}}Q_\nu(\s_{\L}|\eta_\L)\log\frac{\int\nu(d\s) \frac{M_\t(\s_{\L}\xi_{\L^c},\eta_\L)}{M_\t(\s,\eta_\L)}M_\t(\s,\eta_\L)1_{\s_{\L}}(\s)}{\int\mu(d\s) \frac{M_\t(\s_{\L}\xi_{\L^c},\eta_\L)}{M_\t(\s,\eta_\L)}M_\t(\s,\eta_\L)1_{\s_{\L}}(\s)}\frac{M_\t\mu(\eta_\L)}{M_\t\nu(\eta_\L)}\cr
&\leq\sum_{\eta_\L}M_\t\nu(\eta_\L)\sum_{\s_{\L}}Q_\nu(\s_{\L}|\eta_\L)\log\frac{Q_\nu(\s_{\L},\eta_\L)}{M_\t\nu(\eta_\L)}\frac{M_\t\mu(\eta_\L)}{Q_\mu(\s_{\L},\eta_\L)}+\log\frac{\sup_{\s_{\L}=\tilde\s_{\L},\eta}\frac{M_\t(\s,\eta_\L)}{M_\t(\tilde\s,\eta_\L)}}{\inf_{\s_{\L}=\tilde\s_{\L},\eta}\frac{M_\t(\s,\eta_\L)}{M_\t(\tilde\s,\eta_\L)}}\cr
%&=\sum_{\eta_\L}M_\t\nu(\eta_\L)\sum_{\s_{\L}}Q_\nu(\s_{\L}|\eta_\L)\log\frac{Q_\nu(\s_{\L},\eta_\L)}{M_\t\nu(\eta_\L)}\frac{M_\t\mu(\eta_\L)}{Q_\mu(\s_{\L},\eta_\L)}+2\log\sup_{\s_{\L}=\tilde\s_{\L},\eta}\frac{M_\t(\s,\eta_\L)}{M_\t(\tilde\s,\eta_\L)}\cr
%&=\sum_{\eta_\L}M_\t\nu(\eta_\L)\sum_{\s_{\L}}Q_\nu(\s_{\L}|\eta_\L)\log\frac{Q_\nu(\s_{\L}|\eta_\L)}{Q_\mu(\s_{\bar\L}|\eta_\L)}+2\log e^{o(|\L|)}\cr
&=\sum_{\eta_\L,\s_{\L}}Q_\nu(\s_{\L},\eta_\L)\log\frac{ Q_\nu(\s_{\L}|\eta_\L)}{Q_\mu(\s_{\L}|\eta_\L)}+2 o(|\L|)\cr
%&=\sum_{\eta_\L}h_{\L}(Q_\nu(\cdot|\eta_\L)|Q_\mu(\cdot|\eta_\L))M_\t\nu(\eta_\L)+o(|\L|)\cr
\end{split}
\end{equation*}
where we used \eqref{OneStepError} and $\sum_{\eta_\L\in\{1,\dots,q\}^\L}M_\t\nu(\eta_\L)Q_\nu(\s_{\L}|\eta_\L)=\sum_{\eta_\L}Q_\nu(\s_{\L},\eta_\L)=\nu(\s_{\L})$. Notice, we get a similar bound from below, i.e 
\begin{equation*}\label{EntropyProduction}
\begin{split}
&h_{\L}(\nu|\mu)-h_\L(M_\t\nu|\mu)
\geq\sum_{\eta_\L,\s_{\L}}Q_\nu(\s_{\L},\eta_\L)\log\frac{Q_\nu(\s_{\L}|\eta_\L)}{Q_\mu(\s_{\L}|\eta_\L)}-o(|\L|).\cr
\end{split}
\end{equation*} 
%Further notice 
%\begin{equation*}
%\begin{split}
%\int Q_\nu(d\s,d\eta)&\log\frac{\hat Q_\nu(\s_{\L}|\eta_\L)}{\hat Q_\mu(\s_{\L}|\eta_\L)}=\int Q_\nu(d\s,d\eta)Q_\nu(\log\frac{\hat Q_\nu(\s_{\L}|\eta_\L)}{\hat Q_\mu(\s_{\L}|\eta_\L)}|\mathcal{F}_{\L}\times\mathcal{F}_\L)(\s,\eta)\cr
%&=\int Q_\nu(d\s,d\eta)\sum_{\eta_\L,\s_{\L}}Q_\nu(\s_{\L},\eta_\L|\mathcal{F}_{\L}\times\mathcal{F}_\L)(\s,\eta)\log\frac{\hat Q_\nu(\s_{\L}|\eta_\L)}{\hat Q_\mu(\s_{\L}|\eta_\L)}\cr
%&=\sum_{\eta_\L,\s_{\bar\L}}Q_\nu(\s_{\L},\eta_\L)\log\frac{\hat Q_\nu(\s_{\L}|\eta_\L)}{\hat Q_\mu(\s_{\L}|\eta_\L)}\cr
%\end{split}
%\end{equation*}
%and t
Together we have the following identity
\begin{equation}\label{EntropyProduction_2}
\begin{split}
&h_{\L}(\nu|\mu)-h_\L(M_\t\nu|\mu)=\E^{Q_\nu}[\log\frac{Q_\nu(\s_{\L}|\eta_\L)}{Q_\mu(\s_{\L}|\eta_\L)}]\pm o(|\L|)\cr
&=\E^{Q_\nu}[\log\frac{Q_\nu(\s_{\L}|\eta_\L)}{Q_\nu(\s_{\L}|\eta)}]+\E^{Q_\nu}[\log\frac{\hat M_{\t,\nu}(\eta, \s_{\L})}{M_{-\t}(\eta,\s_{\L})}]+\E^{Q_\nu}[\log\frac{Q_\mu(\s_{\L}|\eta)}{Q_\mu(\s_{\L}|\eta_\L)}]\pm o(|\L|).\cr
\end{split}
\end{equation}

Under the volume limit, the l.h.s of \eqref{EntropyProduction_2} becomes the l.h.s of \eqref{Lemma1} and for the second summand on the r.h.s of \eqref{EntropyProduction_2} we have 
\begin{equation*}\label{EntropyProduction_3}
\begin{split}
\limsup_{\L\uparrow G}\frac{1}{|\L|}\E^{Q_\nu}[\log\frac{\hat M_{\t,\nu}(\eta, \s_{\L})}{M_{-\t}(\eta,\s_{\L})}]
%&=\limsup_{\L\uparrow G}\frac{1}{|\L|}\int\nu(d\s)M_\t(\s,d\eta)1_{\s_\L}(\s)\log\frac{\hat M_{\t,\nu}(\eta, \s_{\L})}{M_{-\t}(\eta,\s_{\L})}\cr
%&=\limsup_{\L\uparrow G}\frac{1}{|\L|}\int M_\t\nu(d\eta)\hat M_{\t,\nu}(\eta,d\s)1_{\s_\L}(\s)\log\frac{\hat M_{\t,\nu}(\eta, \s_{\L})}{M_{-\t}(\eta,\s_{\L})}\cr
&=\limsup_{\L\uparrow G}\frac{1}{|\L|}\int M_\t\nu(d\eta)h_\L[\hat M_{\t,\nu}(\eta, \cdot)|M_{-\t}(\eta,\cdot)]\cr
&=\mathcal H_{M_\t\nu}(\hat M_{\t,\nu}|M_{-\t})
\end{split}
\end{equation*}
which is the r.h.s of \eqref{Lemma1}. Hence, in order to prove \eqref{Lemma1} it suffices to show, that the first and the third summand on the r.h.s of \eqref{EntropyProduction_2} are $o(|\L|)$-functions. 
Since the third summand is not a special case of the first summand, we have to proceed separately.
%
%Since for the l.h.s of we have $\limsup_{\L\uparrow\Z^d}[h_{\L}(\nu|\mu)-h_\L(M_\t\nu|\mu)]=h(\nu|\mu)-h(M_\t\nu|\mu)$ which is and 

%\bigskip
%Notice, in \eqref{EntropyProduction_2} the third summand is a special case of the first summand. 

\medskip
\textbf{For the first summand in \eqref{EntropyProduction_2}} we can follow closely the arguments from \cite{DaLoRo02} Proposition 2.1. For the readers convenience we provide them here as well. 
%we can introduce the density w.r.t the (maybe infinite) product of the uniform measure $\l^\L$ on $\{1,\dots,q\}^L$. 
%We denote by $\l^\otimes(\eta_\L)$ the projection of $\l^\otimes$ on $\mathcal F_\L$. 
Let $\{i_1,\dots, i_{|\L|}\}$ be the lexicographic ordering of the elements of $\L$ and define $\L_k=\{i_1,\dots,i_k\}$ for $1\leq k\leq\L$ with $\L_0=\emptyset$. By Bayes' theorem we have 
\begin{equation*}\label{44}
\begin{split}
Q_\nu(\s_{\L}|\eta_\L)
%&=\frac{Q_\nu(\s_{i_1},\eta_\L)}{Q_\nu(\{1,\dots,q\}^G,\eta_\L)}\frac{Q_\nu(\s_{i_1,i_2},\eta_\L)}{Q_\nu(\s_{i_1},\eta_\L)}\cdots\frac{Q_\nu(\s_{\L_{|\L|}},\eta_\L)}{Q_\nu(\s_{\L_{|\L|-1}},\eta_\L)}\cr
&=Q_\nu(\s_{i_1}|\eta_\L)Q_\nu(\s_{i_2}|\eta_\L,\s_{i_1})\cdots Q_\nu(\s_{i_{|\L|}}|\eta_\L,\s_{\L_{|\L|-1}})\cr
\end{split}
\end{equation*}
and hence we can write
\begin{equation*}\label{45}
\begin{split}
\log Q_\nu(\s_{\L}|\eta_\L)=\sum_{k=1}^{|\L|}\log Q_\nu(\s_{i_k}|\eta_\L,\s_{\L_{k-1}}).
\end{split}
\end{equation*}
By translation invariance of $Q_\nu$ we have
\begin{equation*}\label{46}
\begin{split}
\E^{Q_\nu}[\log Q_\nu(\s_{i_k}|\eta_\L,\s_{\L_{k-1}})]=\E^{Q_\nu}[\log Q_\nu(\s_0|\eta_{\theta_{-i_k}\L},\s_{\theta_{-i_k}\L_{k-1}})]
%&=\E^Q\Bigl[\log q\hat Q(\s_0|\eta_{\theta_{-i_k}\L},\s_{\theta_{-i_k}\bar\L_{k-1}})\Bigr]
\end{split}
\end{equation*}
where $\theta_{i}\L$ denotes a lattice translation of the set $\L\subset G$ by $i\in G$. 
We want to use the Shannon-Breiman-McMillan Theorem from \cite{Ba85} as used in \cite{DaLoRo02} and therefore understand the conditional densities as a sequence of densities belonging to a stochastic process under the invariant measure $Q_\nu$.

Let $G_{-}:=\{i\in G:i\prec0\}$, where $\prec$ is the lexicographic order. By the Shannon-Breiman-McMillan Theorem the stationary and ergodic process has the "Asymptotic Equipartition Property" and thus the sequence of conditional measures has an almost sure limit. In particular it is a Cauchy sequence and hence, for every $\e>0$ there are $A\subset G,B\subset G_{-}$ finite such that if $A\subset V$ and $B\subset W\subset G_{-}$ we have
\begin{equation*}\label{Shannon-Breiman-McMillan}
\begin{split}
\Bigr|\E^{Q_\nu}[\log Q_\nu(\s_0|\eta_{V},\s_{W})]-\E^{Q_\nu}[\log Q_\nu(\s_0|\eta_{A},\s_{B})]\Bigl|<\e.
\end{split}
\end{equation*}
Also notice for all $\L$ and $i_k$ there is the simple bound
\begin{equation*}\label{BrutalBound}
\begin{split}
\E^{Q_\nu}[\log Q_\nu(\s_0|\eta_{\theta_{-i_k}\L},\s_{\theta_{-i_k}\L_{k-1}})]\leq\log q.
\end{split}
\end{equation*}
Now we can separate bulk and boundary terms. This gives
\begin{equation}\label{TogetherA}
\begin{split}
%&\limsup_{\L\uparrow G}\frac{1}{|\L|}\E^{Q_\nu}\Bigl[\log\frac{ Q_\nu(\s_{i_k}|\eta_\L,\s_{\L_{k-1}})}{\l(\s_{i_k})}\Bigr]\cr
%&=
&\limsup_{\L\uparrow G}\frac{1}{|\L|}\sum_{k=1}^{|\L|}
\E^{Q_\nu}[\log Q_\nu(\s_0|\eta_{\theta_{-i_k}\L},\s_{\theta_{-i_k}\L_{k-1}})]\cr
&=\limsup_{\L\uparrow G}\frac{1}{|\L|}\Bigl\{\sum_{k:A\subset\theta_{-i_k}\L,B\subset\theta_{-i_k}\L_{k-1}}
\E^{Q_\nu}[\log Q_\nu(\s_0|\eta_{\theta_{-i_k}\L},\s_{\theta_{-i_k}\L_{k-1}})]\cr
&\hspace{5cm}+|\{k:A\not\subset\theta_{-i_k}\L,B\not\subset\theta_{-i_k}\L_{k-1}\}|\log q\Bigr\}.\cr
\end{split}
\end{equation}
For large $\L$ the first term on the r.h.s of \eqref{TogetherA} contains the bulk of the summands. The second term on the r.h.s of \eqref{TogetherA} is of boundary order. Hence
\begin{equation*}\label{Together2}
\begin{split}
&\limsup_{\L\uparrow G}\frac{1}{|\L|}\E^{Q_\nu}[\log Q_\nu(\s_{\L}|\eta_\L)]
\leq\E^{Q_\nu}[\log Q_\nu(\s_0|\eta_{A},\s_{B})]+\e\hspace{1cm}\text{ and }\cr
&\limsup_{\L\uparrow G}\frac{1}{|\L|}\E^{Q_\nu}[\log Q_\nu(\s_{\L}|\eta_\L)]
\geq\E^{Q_\nu}[\log Q_\nu(\s_0|\eta_{A},\s_{B})]-\e.\cr
\end{split}
\end{equation*}
Letting $\e$ go to zero we have 
\begin{equation*}\label{Together1}
\begin{split}
&\limsup_{\L\uparrow G}\frac{1}{|\L|}\E^{Q_\nu}[\log Q_\nu(\s_{\L}|\eta_\L)]
=\E^{Q_\nu}[\log Q_\nu(\s_0|\eta,\s_{G_{-}})].
\end{split}
\end{equation*}
Using the same arguments as above one also shows 
\begin{equation*}\label{Together2}
\begin{split}
&\limsup_{\L\uparrow G}\frac{1}{|\L|}\E^{Q_\nu}[\log Q_\nu(\s_{\L}|\eta)]
=\E^{Q_\nu}[\log Q_\nu(\s_0|\eta,\s_{G_{-}})].
\end{split}
\end{equation*}

\textbf{For the second summand in \eqref{EntropyProduction_2}} recall 
%
%notice first: It is simple to verify that the time reversed transition operator for $M_{\t}$, defined for a starting measure
%% $\nu$ to be 
%%\begin{equation}\label{zwei}
%% \begin{split}
%% &\int \nu(d\eta')\int \mu[\eta'](d\eta)f(T(\eta+\t),\eta')=\int (M_\t\nu)(d\s')\int \hat M(d\eta'|\s') f(\s',\eta')
%%\end{split}
%%\end{equation}
%%is given by $\hat M(d\eta'|\s')=M_{-\t}(d \eta' | \s')$ for all starting measures 
%$\mu_\phi\in\text{ex } \GG_{\theta}(\g')$ is $M_{-\t}$, independently of the angle $\phi$. This is clear since
%%The proof is clear since 
%\begin{equation*}\label{BackwardsGibbs}
% \begin{split}
%\int \mu_\phi(d\eta')\int M_\t(\eta',d\s')f(\s',\eta')&=\int \mu_\phi(d\eta')\int \mu_G[\eta'](d\eta)f(T(\eta+\t),\eta')\cr
%&=\int \mu_\phi(d\eta')\int \mu_G[\eta'](d\eta)f(T(\eta+\t),T(\eta))\cr
%&=\int \mu_\phi(d\eta)f(T(\eta+\t),T(\eta))\cr
%&=\int \mu_{\phi+\t}(d\s)f(T(\s),T(\s-\t))\cr
%&=\int \mu_{\phi+\t}(d\s')\int \mu_G[\s'](d\s)f(\s',T(\s-\t))\cr
%&=\int \mu_{\phi+\t}(d\s')\int M_{-\t}(\s',d\eta') f(\s',\eta').
%\end{split}
%\end{equation*}
%%Further notice by \eqref{BackwardsGibbs} for the backwards transition operator for an extremal translation-invariant Gibbs measure we have 
%%\begin{equation*}\label{BackwardsGibbs2}
%% \begin{split}
%%\hat Q_{\mu}(\s_{\L}|\eta)=M_{-\t}(\eta,\s_{\L}).
%%   \end{split}
%%\end{equation*}
%%where $\hat M$ is a $q$-fold composition of $M_{-\t}$. 
equation \eqref{OneStepError}. Of course also $M_{-\t}$ has such an error bound and we have
% and \eqref{TwoStepError}
\begin{equation*}\label{BackwardsGibbsError}
 \begin{split}
\frac{M_{-\t}(\eta_{\L}\xi_{\L^c},\s_{\L})}{M_{-\t}(\eta,\s_{\L})}\leq e^{o(|\L|)}.
   \end{split}
\end{equation*}
Hence with Lemma \ref{GibbsBackwards} we can write
%Since $\mu$ is assumed to be time-stationary for $M_\t$ and a translation-invariant Gibbs measure, it must be a symmetric convex combination of the extremal shift-invariant Gibbs measures. Hence the time-reversed transition operator for this $\mu$ is also given by $M_{-\t}$.\textbf{ ??????}
\begin{equation*}\label{BackwardsGibbs3}
 \begin{split}
\frac{Q_{\mu}(\s_{\L}|\eta_\L)}{Q_{\mu}(\s_{\L}|\eta)}
%=\frac{\int M\mu(d\tilde\eta|\eta_\L)\hat M(\eta_\L\tilde\eta_{\L^c},\s_\L)}{\hat M(\eta,\s_{\L})}
=\int M_\t\mu(d\tilde\eta|\eta_\L)\frac{M_{-\t}(\eta_\L\tilde\eta_{\L^c},\s_\L)}{M_{-\t}(\eta,\s_{\L})}
\leq e^{o(|\L|)}
   \end{split}
\end{equation*}
and thus $\E^{Q_{\nu}}[\log\frac{Q_{\mu}(\s_{\L}|\eta)}{Q_{\mu}(\s_{\L}|\eta_\L)}]=\pm o(|\L|)$.
As a result, in the infinite-volume limit, the $o(|\L|)$ terms vanish and this completes the proof.
$\Cox$

%\bigskip
%UNKLAR: Wo wird bei DPLR und somit auch bei uns die Gibbs Eigenschaft benutzt?

\bigskip
Recall that our goal is to show that time-stationary measures must be Gibbs measures. The following result is slightly more general. Using the characterization of the entropic loss from the preceding lemma we will show that zero entropic loss implies the single site DLR equation and thus the Gibbs property. The main tools for the proof are an adaptation of the Gibbs variational principle and an argument from \cite{Ku84} to infer the DLR equation from the backwards transition operator. For more general background about the Gibbs variational principle also in the context of generalized Gibbsian measures see \cite{Ge11,KuLeRe04,EnVe04}.

\begin{thm}\label{ZeroEntrGibbs} 
Suppose that $\mu'$ is a translation-invariant Gibbs measure w.r.t the specification $\g'$ and also time-stationary w.r.t $M_\t$ where $0<\t<2\pi/q$. Then, for any translation-invariant measure $\nu'$ with
%the vanishing of the entropic loss under $M_\t$, i.e.  
\begin{equation*}\label{zwei}
 \begin{split}
 &h(\nu'|\mu')=h(M_\t\nu'|\mu' )
   \end{split}
\end{equation*}
we have that $\nu'$ is also a Gibbs measure for the same specification, i.e $\nu'\in \GG_\theta(\g')$.
Conversely if $\nu'\in \GG_\theta(\g')$ then also $h(\nu'|\mu')=h(M_\t\nu'|\mu' )$.
\end{thm}
Please note, Theorem \ref{ZeroEntrGibbs} and the following results are valid for general $\t\in\R$. We set $0<\t<2\pi/q$ only to keep the presentation more transparent.

\bigskip
{\bf Proof: }Let us again suppress all primes in the notation whenever unambiguous.

\textbf{Step 1: }We show that 
\begin{equation}\label{Step1}
\begin{split}
\hat M_{\t,\nu}(\eta,\s_{\L})=M_{-\t}(\eta,\s_{\L})\hspace{0.5cm}\text{for }M_\t\nu\text{-a.a. }\eta
\end{split}
\end{equation}
under the hypothesis of the theorem. To do so we extend the proof of the Gibbs variational principle in \cite{Ge11} Chapter 15.4. to our situation where we have to deal with an additional dependence on the future configuration $\eta$ distributed according to $M_\t\nu$.
Define $h_\L^\eta(\hat M_\nu|\hat M_\mu):=h_\L[\hat M_{\t,\nu}(\eta,\cdot)|M_{-\t}(\eta,\cdot)]$, then by Fatou's Lemma and Lemma \ref{Extension of Proposition 2.1. from DaPrLoRo02} we have
\begin{equation}\label{Zero}
 \begin{split}
0=\HH_{M_\t\nu}(\hat M_\nu|\hat M_\mu)&=\int \nu(d\eta) \limsup_{\L}\frac{1}{|\L|}h_\L^\eta(\hat M_\nu|\hat M_\mu)\cr
&\geq\limsup_{\L}\frac{1}{|\L|}\int M_\t\nu(d\eta) h_\L^\eta(\hat M_\nu|\hat M_\mu)=0.\cr
%&\geq\liminf_{\L}\frac{1}{|\L|}\int M_\t\nu(d\eta) h_\L^\eta(\hat M_\nu|\hat M_\mu)=0.\cr
   \end{split}
\end{equation}
In particular for all cofinal sequences of finite sets
$0=\lim_{\L_n\uparrow G}\frac{1}{|\L_n|}\int M_\t\nu(d\eta) h_{\L_n}^\eta(\hat M_\nu|\hat M_\mu)$. 
%First we prove that $Q_\nu$-a.s. for all finite sets $\L$ we have $Q_\nu(\s_{\L}|\s_{\L^c},\eta)=Q_\mu(\s_{\L}|\s_{\L^c},\eta)$.
Since the local relative entropy $h_\L^\eta(\hat M_\nu|\hat M_\mu)$ must be finite for all $\L$ and $M_\t\nu$-a.a. $\eta$ there exists a density 
 \begin{equation*}\label{vier2}
 \begin{split}
f_\L^\eta:=\frac{d\hat M_\nu(\eta, \cdot |_{\L})}{d\hat M_\mu(\eta, \cdot |_{\L})}\cr
   \end{split}
\end{equation*}
depending on configurations inside $\L$ only.

\medskip
\textbf{Step 1a: }We show that expectations of local entropy densities behave nice w.r.t volumes in the following sense:
For any $\d>0$ and any cube $C\supset\L$ there exists a finite set $\D$ with 
$\D\supset C$ such that 
\begin{equation*}\label{funf}
 \begin{split}
\int M_\t\nu(d\eta) h_{\D}^\eta(\hat M_\nu|\hat M_\mu)-\int M_\t\nu(d\eta) h_{\D\setminus\L}^\eta(\hat M_\nu|\hat M_\mu)\leq\d.\cr
   \end{split}
\end{equation*}
The proof is an integrated version of the first step of the proof of the variational principle in \cite{Ge11} Chapter 15.4. Indeed, by \eqref{Zero} we can find $n\geq1$ such that for the centered hypercube $\L_n$ we have $|\L_n|\geq|C|$ and $|\L_n|^{-1}\int M_\t\nu(d\eta) h_{\L_n}^\eta(\hat M_\nu|\hat M_\mu)\leq  \d/2^d|C|$. Now we can choose an integer $m\geq1$ in such a way that
\begin{equation*}\label{sechs}
 \begin{split}
m^d|C|\leq |\L_n|\leq (2m)^d|C|.
   \end{split}
\end{equation*}
Further let us choose $m^d$ lattice sites $i(1),\dots, i(m^d)$ in such a way that the translates $C(k) = C + i(k)$, $1\leq k\leq m^d$ are pairwise disjoint subsets of $\L_n$. For 
each $1\leq k\leq m^d$ we put $W(k) = C(1)\cup\dots\cup C(k)$ and $\L(k) = \L + i(k)$. Then 
using the monotonicity of the relative density we can write 
\begin{equation*}\label{sieben}
 \begin{split}
\frac{1}{m^d}&\sum_{i=1}^{m^d}[\int M_\t\nu(d\eta) h_{W(k)}^\eta(\hat M_\nu|\hat M_\mu)-\int M_\t\nu(d\eta) h_{W(k)\setminus\L(k)}^\eta(\hat M_\nu|\hat M_\mu)]\cr
&\leq\frac{1}{m^d}\sum_{i=1}^{m^d}[\int M_\t\nu(d\eta) h_{W(k)}^\eta(\hat M_\nu|\hat M_\mu)-\int M_\t\nu(d\eta) h_{W(k)\setminus C(k)}^\eta(\hat M_\nu|\hat M_\mu)]\cr
&=\frac{1}{m^d}\int M_\t\nu(d\eta) h_{W(m^d)}^\eta(\hat M_\nu|\hat M_\mu)\cr
&\leq\frac{1}{m^d}\int M_\t\nu(d\eta) h_{\L_n}^\eta(\hat M_\nu|\hat M_\mu)\leq 2^d|C||\L_n|^{-1}\int M_\t\nu(d\eta) h_{\L_n}^\eta(\hat M_\nu|\hat M_\mu)\leq\d.\cr
   \end{split}
\end{equation*}
Consequently, there exists an index $k$ such that $$\int M_\t\nu(d\eta) h_{W(k)}^\eta(\hat M_\nu|\hat M_\mu)-\int M_\t\nu(d\eta) h_{W(k)\setminus\L(k)}^\eta(\hat M_\nu|\hat M_\mu)\leq\d.$$ The claim of Step 1a thus follows by putting $\D:= W(k)-i(k)$ and using the translation-invariance of $\hat M_\nu(\eta, \cdot |_{\L})$ and $\hat M_\mu(\eta, \cdot |_{\L})$ under $M_\t\nu(d\eta)$. 

\medskip
\textbf{Step 1b: }We want to show that closeness of integrated entropy densities implies closeness of integrated densities. This again is an integrated version of the analogous statement in \cite{Ge11} Chapter 15.4. More precisely we show that for each $\e>0$ there exists some $\d>0$ such that $$\int M_\t\nu(d\eta)\int\hat M_\mu(d\s,\eta)|f_\D^\eta(\s)-f_{\D\setminus\L}^\eta(\s)|\leq\e$$ whenever $\L\subset\D$ and $$\int M_\t\nu(d\eta) h_{\D}^\eta(\hat M_\nu|\hat M_\mu)-\int M_\t\nu(d\eta) h_{\D\setminus\L}^\eta(\hat M_\nu|\hat M_\mu)\leq\d.$$
Indeed, let $\psi(x):=x\log x-x+1$ for $x\geq0$, then
\begin{equation*}\label{acht}
 \begin{split}
&\int M_\t\nu(d\eta) h_{\D}^\eta(\hat M_\nu|\hat M_\mu)-\int M_\t\nu(d\eta) h_{\D\setminus\L}^\eta(\hat M_\nu|\hat M_\mu)\cr
&=\int M_\t\nu(d\eta)\int \hat M_\nu(\eta,d\s)\log\frac{f_\D^\eta(\s)}{f_{\D\setminus\L}^\eta(\s)}\cr
&=\int M_\t\nu(d\eta)\int \hat M_\mu(\eta,d\s)f_{\D\setminus\L}^\eta(\s)\psi(\frac{f_\D^\eta(\s)}{f_{\D\setminus\L}^\eta(\s)}).\cr
   \end{split}
\end{equation*}
Notice, there is a number $0<r<\infty$ such that $|x-1|\leq r\psi(x) + \e/2$ for all $x\geq0$. By putting $\d = \e/2r$ we can thus write
\begin{equation*}\label{neun}
 \begin{split}
&\int M_\t\nu(d\eta)\int\hat M_\mu(d\s,\eta)|f_\D^\eta(\s)-f_{\D\setminus\L}^\eta(\s)|\cr
&=\int M_\t\nu(d\eta)\int\hat M_\mu(d\s,\eta)f_{\D\setminus\L}^\eta(\s)|\frac{f_\D^\eta(\s)}{f_{\D\setminus\L}^\eta(\s)}-1|\cr
&\leq r\int M_\t\nu(d\eta)\int \hat M_\mu(\eta,d\s)f_{\D\setminus\L}^\eta(\s)\psi(\frac{f_\D^\eta(\s)}{f_{\D\setminus\L}^\eta(\s)})+\frac{\e}{2}\leq\e.\cr
   \end{split}
\end{equation*}

\medskip
\textbf{Step 1c: }To prove the desired result for the Step 1, i.e that $M_\t\nu$-a.s. we have $\hat M_{\t,\nu}(\eta,\s_{\L})=M_{-\t}(\eta,\s_{\L})$, 
%In order to conclude the desired identity \eqref{Step1} 
we show the DLR equation for the backwards operator corresponding to $\nu$ for the specification given by the backwards transition operator corresponding to the Gibbs measure $\mu.$ Indeed, let $g$ be a local test-function and $\e>0$. 
%Notice, since $\mu$ is a translation-invariant and time-stationary Gibbs measure, $\hat M_{\t,\mu}(\eta, d\s)=M_{-\t}(\eta, d\s)$ 
%where 
%\begin{equation}\label{hatM}
% \begin{split}
%\hat M(\s',\eta'_{\L})&=\int M_{-\t}(\s',d\eta'_1)\int M_{-\t}(\eta'_1,d\eta'_2)\cdots\int M_{-\t}(\eta'_{q-1},d\eta'_q) M_{-\t}(\eta'_q,d\eta'_\L)\cr
%  \end{split}
%\end{equation}
%independent of $\mu$. 
Since $\mu$ is a Gibbs measure 
%and $M_{\t}$ is assumed to be quasilocal, 
there exists a local function $\tilde g^\eta$ with support independent of $\eta$ (which depend only on configurations outside $\L$) such that 
$\sup_{\eta}\Vert\g^\eta_\L(g|\cdot)-\tilde g^\eta(\cdot)\Vert<\e$ where $$\g^\eta_\L(d\tilde\s|\s_{\L^c}):=Q_\mu(d\tilde\s|\s_{\L^c},\eta).$$
To see this we note that $M_{-\t}(\eta,d\s)$ as a measure on the $\s$ is in the Dobrushin region uniformly in $\eta$. 
This is true by arguments provided in the proof of Step 2 below equation \eqref{TauCoarse}.
Let $C\supset \L$ be a cube such that $g$ depends only on configuration on $C$ and $\tilde g^\eta$ depends only on $C\setminus\L$. Choose $\d$ in terms of $\e$ as in Step 1b, and define $\D$ in terms of $C$ and $\d$ as in Step 1a. 
%Then we can estimate the difference between the left and right side of the DLR equation in the following way: 
\begin{equation*}\label{zehn}
 \begin{split}
&\int M_\t\nu(d\eta)|\int\hat M_{\t,\nu}(\eta,d\s)g(\s)-\int\hat M_{\t,\nu}(\eta,d\s)\g^\eta_\L(g|\s)|\cr
&\leq\int M_\t\nu(d\eta)\Bigl[\int\hat M_{\t,\nu}(\eta,d\s)|\g^\eta_\L(g|\s)-\tilde g^\eta(\s)|\cr
&\hspace{1cm}+|\int\hat M_{\t,\nu}(\eta,d\s)\tilde g^\eta(\s)-\int M_{-\t}(\eta,d\s)f^\eta_{\D\setminus\L}(\s)\tilde g^\eta(\s)|\cr
&\hspace{1cm}+\int M_{-\t}(\eta,d\s)f^\eta_{\D\setminus\L}(\s)|\tilde g^\eta(\s)-\g^\eta_\L(g|\s)|\cr
&\hspace{1cm}+|\int M_{-\t}(\eta,d\s)(f^\eta_{\D\setminus\L}(\s)(\g^\eta_\L(g|\s)-g(\s))|\cr
&\hspace{1cm}+\Vert g\Vert\int M_{-\t}(\eta,d\s)|f^\eta_{\D\setminus\L}(\s)-f^\eta_{\D}(\s)|\cr
&\hspace{1cm}+|\int M_{-\t}(\eta,d\s)f^\eta_{\D}(\s)g(\s)-\int\hat M_{\t,\nu}(\eta,d\s)g(\s)|\Bigr]\cr
   \end{split}
\end{equation*}
Since $\tilde g^\eta$ depends only on $\D\setminus\L$ and $g$ depends only $\D$, the second and the last term on the right are zero. The fourth term vanishes because $\int M_{-\t}(\eta,d\s)\int\g^\eta(d\tilde\s|\s)=\int M_{-\t}(\eta,d\s)$ and $f^\eta_{\D\setminus\L}$ depends only on $\L^c$. Due to the choice of $g$, the first and the third term are each at most $\e$. The only non-trivial term is the fifth one. This term is not larger than $\Vert g\Vert\e$ because of 
our choice of $\D$ see Step 1b. As $\e$ was arbitrary, we conclude that $\int M_\t\nu(d\eta)|\int\hat M_{\t,\nu}(\eta,d\s)g(\s)-\int\hat M_{\t,\nu}(\eta,d\s)\g^\eta_\L(g|\s)|=0$. Hence $M_\t\nu$-a.s. we have $\int\hat M_{\t,\nu}(d\s,\eta)g(\s)=\int\hat M_{\t,\nu}(d\s,\eta)\g^\eta_\L(g|\s)$. And thus $\hat M_{\t,\nu}$ is the unique Gibbs measure for $\g_\L^\eta$.
%Setting $g=1_{\xi_\L}$ and $\bar\L$ large we have $M_\t\nu$-a.s. $\hat M_{\t,\nu}(\eta,\xi_\L)=\int\hat M_{\t,\nu}(d\s,\eta)Q_\mu(1_{\xi_\L}|\s_{\bar\L^c},\eta)$. Letting $\bar\L\uparrow G$ and using again the quasilocality and the uniqueness of the Gibbs measure corresponding to $\g_\L^\eta$, 
In particular equation \eqref{Step1} follows.

\bigskip
\textbf{Step 2: }
As a consequence of Step 1, equation \eqref{Step1} holds and we can conclude that $Q_\nu$-a.s.
\begin{equation*}\label{Together_5}
 \begin{split}
Q_\nu(\s_{\L}|\s_{\L^c},\eta)=Q_\mu(\s_{\L}|\s_{\L^c},\eta).
   \end{split}
\end{equation*}
In what follows, we use the basic idea from \cite{Ku84} to infer from this the single-site DLR equation for $\nu$. Notice, for $i\in G$, finite sets $\L\subset G\setminus i$, $\bar\L\subset G$ and $\tilde\s_{i^c}=\s_{i^c}$ we have
\begin{equation*}\label{Together_6}
\begin{split}
\mu(\s_i|\s_{\L})Q_\mu(\eta_{\bar\L}|\s_{\L\cup i})Q_\mu(\tilde\s_{i}|\s_{\L},\eta_{\bar\L})=
\mu(\tilde\s_i|\s_{\L})Q_\mu(\eta_{\bar\L}|\tilde\s_{\L\cup i})Q_\mu(\s_{i}|\s_{\L},\eta_{\bar\L})
\end{split}
\end{equation*}
where $\mu$ is a Gibbs measure for the specification $\g'$ given in \eqref{Coarse_Specification}. 
%given by
%\begin{equation*}\label{Specification}
% \begin{split}
%\g'_\L(\o_\L'|\o'_{\L^c})=\frac{\mu[\o'_{\L^c}](\l^{\L}(e^{-H_\L}1_{\s'_{\L}}))}{\mu[\o'_{\L^c}](\l^{\L}(e^{-H_\L}))}.
%   \end{split}
%\end{equation*}
%Recall that $\g'$ is positive and quasilocal since $\s'\mapsto\mu_G[\s']$ is continuous w.r.t the product-topology and the weak-topology (for more details see \cite{JaKu12} Section 2). 
Letting $\L$ go to $G\setminus i$ and $\bar\L$ go to $G$ we get $Q_\mu$-a.s.
%\begin{equation}\label{Together_7}
%\begin{split}
%\frac{\g'(\s'_i|\s'_{i^c})}{\g'(\tilde\s'_i|\s'_{i^c})}=\lim_{\bar\L\uparrow\Z^d}\frac{M(\tilde\s',\eta_{\bar\L}')}{M(\s',\eta_{\bar\L}')}\frac{\hat Q_\mu(\s'_{i}|\s'_{i^c},\eta')}{\hat Q_\mu(\tilde\s'_{i}|\s'_{i^c},\eta')}
%\end{split}
%\end{equation}
\begin{equation}\label{Kuensch1}
\begin{split}
\frac{\mu(\s_i|\s_{i^c})}{\mu(\tilde\s_i|\s_{i^c})}=\frac{\g'(\s_i|\s_{i^c})}{\g'(\tilde\s_i|\s_{i^c})}=\lim_{\bar\L\uparrow G}\frac{M_\t(\tilde\s,\eta_{\bar\L})}{M_\t(\s,\eta_{\bar\L})}\frac{Q_\mu(\s_{i}|\s_{i^c},\eta)}{Q_\mu(\tilde\s_{i}|\s_{i^c},\eta)}
\end{split}
\end{equation}
where we note that $Q_\mu$-a.s. $\eta$ is accessible for $\s,\tilde\s$ and vice versa. Also we will prove right below that the 
%all denominators are non-zero 
$\lim_{\bar\L\uparrow G}\frac{M_\t(\tilde\s,\eta_{\bar\L})}{M_\t(\s,\eta_{\bar\L})}$ exists
% $Q_\mu$-a.s.
. In the same way one gets $Q_\nu$-a.s.
\begin{equation}\label{Kuensch2}
\begin{split}
\frac{\nu(\s_i|\s_{i^c})}{\nu(\tilde\s_i|\s_{i^c})}=\lim_{\bar\L\uparrow G}\frac{M_\t(\tilde\s,\eta_{\bar\L})}{M_\t(\s,\eta_{\bar\L})}\frac{Q_\nu(\s_{i}|\s_{i^c},\eta)}{Q_\nu(\tilde\s_{i}|\s_{i^c},\eta)}
\end{split}
\end{equation}
and hence $Q_\nu$-a.s. $\frac{\nu(\s_i|\s_{i^c})}{\nu(\tilde\s_i|\s_{i^c})}=\frac{\g'(\s_i|\s_{i^c})}{\g'(\tilde\s_i|\s_{i^c})}$. Summing over the $\tilde\s_i$ we get the desired single-site DLR-equation $$\nu(\s_i|\s_{i^c})=\g'(\s_i|\s_{i^c})$$ $Q_\nu$-a.s.  and thus $\nu$ is a Gibbs measure for $\g'$. 

\medskip
%In case of a local PCA the existence of the limit is simple. 
What remains to be proved is that the limit in \eqref{Kuensch1} and \eqref{Kuensch2} indeed exists. Let us therefore consider $M_\t$, $\tilde\s'_{i^c}=\s'_{i^c}$ and $\eta'_{\bar\L}$ such that $M_\t(\s',\eta'_{\bar\L})\neq0\neq M_\t(\tilde\s',\eta'_{\bar\L})$ for all $\bar\L$. Then we have, using the notation from \eqref{OneStepUpdate}
\begin{equation*}\label{Together_9}
\begin{split}
\frac{M_\t(\tilde\s',\eta_{\bar\L}')}{M_\t(\s',\eta_{\bar\L}')}&=\frac{\mu_{G\ba i}[\s'_{G\ba i}](\l(e^{-H_i}1_{\tilde\s'_i}1_{\eta'_{\bar\L,\t}}))}{\mu_{G\ba i}[\s'_{G\ba i}](\l(e^{-H_i}1_{\s'_i}1_{\eta'_{\bar\L,\t}}))}\frac{\mu_{G\ba i}[\s'_{G\ba i}](\l(e^{-H_i}1_{\s'_i}))}{\mu_{G\ba i}[\s'_{G\ba i}](\l(e^{-H_i}1_{\tilde\s'_i}))}\cr
&=\frac{\mu_{G\ba i}[\s'_{G\ba i}](\l(e^{-H_i}1_{\tilde\s'_i}1_{\eta'_{i,\t}})|1_{\eta'_{\bar\L\setminus i,\t}})}{\mu_{G\ba i}[\s'_{G\ba i}](\l(e^{-H_i}1_{\s'_i}1_{\eta'_{i,\t}})|1_{\eta'_{\bar\L\setminus i,\t}})}\frac{\g'(\s'_i|\s'_{G\ba i})}{\g'(\tilde\s'_i|\s'_{G\ba i})}.\cr
\end{split}
\end{equation*}
Now since we are in the Dobrushin region already for the coarse-graining $T$ with $q\geq q_0(\Phi)$, to coarse-grain with 
\begin{equation}\label{TauCoarse}
\begin{split}
[0,2\pi)\mapsto \{[0,\t),[\t,2\pi/q),[2\pi/q,2\pi/q+\t),\dots,[2\pi-\t,2\pi)\}
\end{split}
\end{equation}
is even finer and thus the conditional first-layer specifications are even more in the Dobrushin region. The idea here is that by conditioning to configurations of small segments of the sphere, one can control the possible interaction strength between spins in such a way, that the conditional specification in the first-layer model is in a parameter regime where there is a unique conditional Gibbs measure. Since the individual segments of the sphere under equal-arc discretization with $q$ become even smaller when we discretize the sphere as in \eqref{TauCoarse} a fortiori the Dobrushin condition of weak interaction is satisfied. In particular there is a unique limiting measure on the conditioning $\s'_{G\ba i}\cap\eta'_{{G\ba i},\t}$. In other words
\begin{equation*}\label{Together_10}
\begin{split}
\lim_{\bar\L\uparrow G}&\frac{\mu_{G\ba i}[\s'_{G\ba i}](\l(e^{-H_i}1_{\tilde\s'_i}1_{\eta'_{i,\t}})|1_{\eta'_{\bar\L\setminus i,\t}})}{\mu_{G\ba i}[\s'_{G\ba i}](\l(e^{-H_i}1_{\s'_i}1_{\eta'_{i,\t}})|1_{\eta'_{\bar\L\setminus i,\t}})}\cr
&=\frac{\mu_{G\ba i}[\s'_{G\ba i}\cap\eta'_{{G\ba i},\t}](\l(e^{-H_i}1_{\tilde\s'_i}1_{\eta'_{i,\t}}))}{\mu_{G\ba i}[\s'_{G\ba i}\cap\eta'_{{G\ba i},\t}](\l(e^{-H_i}1_{\s'_i}1_{\eta'_{i,\t}}))}
%=:F_{\tilde\s'_i}(\s',\eta')
\end{split}
\end{equation*}
exists where the limiting measure $\mu_{G\ba i}[\s'_{G\ba i}\cap\eta'_{{G\ba i},\t}]$ is the unique Gibbs measure conditional to the configuration $\s'_{G\ba i}\cap\eta'_{{G\ba i},\t}$ which is in the $\t$-dependent coarse-graining \eqref{TauCoarse}. 

\medskip
\textbf{Step 3: }The converse statement follows from the fact that if $\nu'\in \GG_\theta(\g')$ then $\hat M_{\t,\nu'}=M_{-\t}$ by Lemma \ref{GibbsBackwards} and hence the l.h.s of equation \eqref{Lemma1} in Lemma \ref{Extension of Proposition 2.1. from DaPrLoRo02} is zero. We conclude $h(\nu'|\mu')=h(M_\t\nu'|\mu' )$.
$\Cox$ 

\bigskip
Notice, the uniform mixture $\mu'_*$ as in \eqref{MixedMeasure} is a Gibbs measure for $\g'$, translation-invariant and time-stationary. Hence we can apply Theorem \ref{ZeroEntrGibbs} to our model with $\mu'=\mu'_*$ and the following corollary holds.

%translation-invariant and time-stationary measures must be Gibbs measures for $\g'$ as well. 

\begin{cor}\label{MainCorollary}
Take $0<\t<\frac{2\pi}{q}$ fixed,
%, maybe sufficiently small. 
%Let $M$ be the $q$-fold iterate of $M_\t$ (to make it non-degenerate and have locally 
%all transitions possible). 
then the lattice-translation invariant measures which are invariant under the one time-step updating with the transition operator $M_\t$ must be contained in the lattice-translation invariant measures $\GG_{\th}(\g')$.
\end{cor} 
Let us come to the main conclusion regarding ergodicity and uniqueness of time-stationary measures for our model.

\section{Stationary measures and rotating states}
Notice that we have assembled essentially two properties of our process. In simple words those are: Translation-invariant time-stationary measures must be Gibbs measures and translation-invariant Gibbs measures rotate. The rotation property alone already induces non-ergodicity. Using the combination of the two properties we can draw conclusions about the number of time-stationary measures and  prove part 3 in Theorem \ref{PCA_Theorem}. 
%The implications of the two properties combined are the content 
%The following theorem makes precises the implications on
%use our knowledge of the behavior of the process on the set of translation-invariant Gibbs measures $\GG_{\th}(\g')$.
\begin{thm}\label{IrrationRational}
Let $0<\t<\frac{2\pi}{q}$ then there are two scenarios.
\begin{enumerate}
\item If $\frac{\t}{2\pi}$ is rational, the Markov process with transition operator $M_\t$ has a continuum of translation-invariant and time-stationary measures. Further the dynamics is non-ergodic with periodic orbits given by $$\O_\a:=\{\mu'_\phi\in\text{ex }\GG_{\th}(\g'):\phi=n\t+\a\text{ for some }n\in\N\}.$$
\item If $\frac{\t}{2\pi}$ is irrational, the process has a unique translation-invariant and time-stationary measure $\mu_*'=\frac{1}{2\pi}\int_0^{2\pi}d\phi\mu'_\phi$ where $\mu'_\phi\in\text{ex }\GG_{\th}(\g')$. The dynamics is non-ergodic since $M^n_\t \mu'_{\phi} \not\rightarrow \mu'_*$ as $n\to\infty$. 
All orbits are dense in %Further any orbit is weakly dense in 
$\text{ex }\GG_{\th}(\g')$ w.r.t the weak topology. 
%is a dense orbit in the sense that for all $0\leq\phi,\psi<2\pi$ there exists a subsequence $(n_k)_{k\in\N}$ such that $M^{n_k}_\t \mu'_{\psi}\to\mu'_{\phi}$ weakly.
%
%, but dense orbits $k\mapsto \mu'_{\phi + k \t q}$. In particular it is also non-ergodic in the sense that  with discrete time $k$ tending to infinity.
\end{enumerate}
\end{thm}
\textbf{Proof: }By Corollary \ref{MainCorollary} the investigation of translation-invariant and time-stationary measures can be restricted to the set $\GG_{\th}(\g')$.

In case $\frac{\t}{2\pi}$ is rational there exist integers $n,m$ such that $n\t=m2\pi$. Hence for all $0\leq\a<2\pi$ the equal weight measure $\frac{1}{n}\sum_{k=1}^{n}\mu'_{k\t+\a}$
%\begin{equation}\label{ShiftedMixedMeasures}
%\begin{split}
%\end{split}
%\end{equation}
is time-stationary where $\mu'_{k\t+\a}$ denotes the extremal translation-invariant Gibbs measure with order parameter $(k\t+\a)\text{mod}_{2\pi}\in S^1$. Further we have $M_\t^{n}\mu'_\phi=\mu'_{\phi+n\t}=\mu'_{\phi+m2\pi}=\mu'_\phi$ for all $\mu'_\phi\in\text{ex }\GG_{\th}(\g')$ and hence for a given starting measure $\mu'_\phi$ the set $\{\mu'_\psi\in\text{ex }\GG_{\th}(\g'): \psi=n\t\phi\text{ for some }n\in\N\}$ is a periodic orbit for the dynamics. In particular the process can not be ergodic.

\medskip
If $\frac{\t}{2\pi}$ is irrational the only possible symmetrically mixed measure as in the first case
%\eqref{ShiftedMixedMeasures} 
is the measure $\mu'_*$. A more detailed measure-theoretic proof for this fact can be found in \cite{JaKu12} Proposition 5.1 where one can replace the Markov semigroup by $(M_\t^n)_{n\in\N}$. For the second statement notice, $M^n_\t \mu'_{\phi}=\mu'_{\phi+n\t}\in\text{ex }\GG_{\th}(\g')$ for all $n\in\N$, but $\mu'_*\not\in\text{ex }\GG_{\th}(\g')$ and hence $M^n_\t \mu'_{\phi} \not\rightarrow \mu'_*$. For the third statement realize that for any $0\leq\phi,\psi<2\pi$ there exists a subsequence such that $(\psi+n_k\t)\text{mod}_{2\pi}\to\phi$ for $k\to\infty$. Hence it suffices to show $\mu'_\phi\to\mu'_0$ weakly for $\phi\to0$. But this is true since for the test-functions $1_{\o'_\L}$ we have 
\begin{equation*}\label{WeakConv}
\begin{split}
\lim_{\phi\to0}\mu'_\phi(\o'_\L)&=\lim_{\phi\to0}\int\mu_\phi(d\o)1_{\o'_\L}(\o)
%=\lim_{\phi\to0}\int\mu_0(d\o)1_{\o'_\L}(\o+\phi)\cr
=\int\mu_0(d\o)\lim_{\phi\to0}\g_{\L}(1_{\o'_\L-\phi}|\o_{\L^c})\cr
%&=\int\mu_0(d\o)\lim_{\phi\to0}\frac{\int_{\o'_\L|^l+\phi}^{\o'_\L|^r+\phi}e^{-H_\L(\s_\L\o_{\L^c})}d\s_\L}{\int e^{-H_\L(\s_\L\o_{\L^c})}d\s_\L}
\end{split}
\end{equation*}
where we used the DLR-equation and the dominated convergence theorem. Now by the continuity of the Hamiltonian and the fact, that the a priori measures is the Lebesgue measure on $S^1$ we have $\lim_{\phi\to0}\g_{\L}(1_{\o'_\L-\phi}|\o_{\L^c})=\g_{\L}(1_{\o'_\L}|\o_{\L^c})$ for any boundary condition $\o_{\L^c}$. In particular the process can again not be ergodic.
$\Cox$

\bigskip
\textbf{Proof of Theorem \ref{PCA_Theorem}: }
The proofs of the locality properties in the past and in the future are given in Propositions \ref{Mixing_Past} and \ref{Mixing_Futur}. For $0\leq\t\leq2\pi/q$ the updating is Bernoulli as can be seen from the representation \eqref{rewriting}. Further defining $\II:=\text{ex }\GG_{\th}(\g')$ in Corollary \ref{QuasilocalBijection} we prove the desired properties of the function $\psi$ (this is part 1). The rotation property is proved to hold in Proposition \ref{RotationProp} (this is part 2). The two possible scenarios for time-stationary measures and periodic orbits are proved in Theorem \ref{IrrationRational} (this is part 3).
$\Cox$

\end{document}